\documentclass[12pt]{article}
\usepackage{amsmath}
\usepackage{amssymb}
\usepackage{xypic}
\begin{document}

\newfont{\matha}{msbm10 scaled \magstep1}
\newfont{\mathb}{msbm8}
\newfont{\mathc}{eufm10 scaled \magstep1}
\newcommand{\Cbold}{\mbox{{\matha C}}}
\newcommand{\Fbold}{\mbox{{\matha F}}}
\newcommand{\Nbold}{\mbox{{\matha N}}}
\newcommand{\Qbold}{\mbox{{\matha Q}}}
\newcommand{\Rbold}{\mbox{{\matha R}}}
\newcommand{\Zbold}{\mbox{{\matha Z}}}

\newcommand{\smNbold}{\mbox{{\mathb N}}}

\newcommand{\gbold}{\mbox{{\mathc g}}}
\newcommand{\nbold}{\mbox{{\mathc n}}}
\newcommand{\slalg}{\mbox{{\mathc sl}}}
\newcommand{\abold}{\mbox{\textbf{a}}}
\newcommand{\bbold}{\mbox{\textbf{b}}}
\newcommand{\cbold}{\mbox{\textbf{c}}}
\newcommand{\dbold}{\mbox{\textbf{d}}}
\newcommand{\ebold}{\mbox{\textbf{e}}}
\newcommand{\hbold}{\mbox{\textbf{h}}}
\newcommand{\ibold}{\mbox{\textbf{i}}}
\newcommand{\kbold}{\mbox{\textbf{k}}}
\newcommand{\lbold}{\mbox{\textbf{l}}}
\newcommand{\mbold}{\mbox{\textbf{m}}}
\newcommand{\pbold}{\mbox{\textbf{p}}}
\newcommand{\tbold}{\mbox{\textbf{t}}}
\newcommand{\ubold}{\mbox{\textbf{u}}}
\newcommand{\vbold}{\mbox{\textbf{v}}}
\newcommand{\wbold}{\mbox{\textbf{w}}}
\newcommand{\xbold}{\mbox{\textbf{x}}}

\newcommand{\gothicS}{\mbox{\mathc{S}}}

\newcommand{\smcbold}{\mbox{\small \textbf{c}}}
\newcommand{\smibold}{\mbox{\small \textbf{i}}}
\newcommand{\smtbold}{\mbox{\small \textbf{t}}}
\newcommand{\smwbold}{\mbox{\small \textbf{w}}}

\newcommand{\kQmod}{\mbox{$kQ$-mod}}
\newcommand{\Etilde}{\mbox{$\tilde{E}$}}
\newcommand{\Ftilde}{\mbox{$\tilde{F}$}}

\title{On crystal operators in Lusztig's parametrizations and 
string cone defining inequalities}
\author{Shmuel \textsc{ZELIKSON}
}
\date{}
\maketitle

\begin{center}
\begin{minipage}{0.9\linewidth}
{\small \textsc{Abstract :} 
Let $\mathbf{w}_0$ be a reduced expression for
the longest element of the Weyl group, adapted
to a quiver of type $A_n$. 
We compare Lusztig's and \linebreak[4] Kashiwara's (string)
parametrizations of the canonical basis 
associated with $\mathbf{w}_0$.
Crystal operators act in a finite number of 
patterns in Lusztig's parametrization, 
which may be seen as vectors.  
We show this set gives  the system of defining \linebreak[4] inequalities
of the string cone constructed by Gleizer and Postnikov. 
We use \linebreak[4] combinatorics of
Auslander-Reiten quivers, and as a by-product
we get an alternative enumeration of 
a set of inequalities defining the string cone, 
based on hammocks.}
\end{minipage}
\end{center}

\section*{1. Introduction}

\hspace*{5mm} Let $U_q(\gbold)$ be the quantized enveloping algebra 
corresponding to a Dynkin diagram
$D$ of type $A_n, D_n, E_n$, defined over $\Cbold(q), \, q$ being an indeterminate. There
 is a braid group action on $U_q(\gbold)$ which enables to construct PBW bases \cite{LusztigBase}
of the positive part $U_q(\nbold^+)$ of $U_q(\gbold)$ . Such a basis $\mathcal{P}_{\mathbf{w}_0}$ depends on the choice of a reduced expression $\mathbf{w}_0$ of the longest element $w_0$ of 
the Weyl group $W$ associated with $D$. However, it was observed by Lusztig 
that the $\Zbold[q^{-1}]$-module $\mathcal{L}$ generated by  $\mathcal{P}_{\mathbf{w}_0}$
is independent of $\mathbf{w}_0$. Furthermore, the image of $\mathcal{P}_{\mathbf{w}_0}$
under the projection $\pi \; : \; \mathcal{L} \longrightarrow \mathcal{L}/q^{-1} \mathcal{L}$ 
is a $\Zbold$-basis $B$ of $\mathcal{L}/q^{-1} \mathcal{L}$, which is again independent of $\mathbf{w}_0$. 
There is a unique basis of $\mathcal{L}$, which is invariant under the $\Cbold$-algebra involution
of $U_q(\nbold^+)$ preserving the generators of $U_q(\nbold^+)$, and sending 
$q$ to $q^{-1}$, and whose image under $\pi$ is $B$. This is the canonical 
basis $\mathcal{B}_{can}$ of Lusztig and Kashiwara \cite{LusztigBase}, \cite{KashiwaraBase}.
This basis is in one-to-one correspondence with any PBW basis $\mathcal{P}_{\mathbf{w}_0}$,
yet  is independent of the choice of $\mathbf{w}_0$. It has with many
remarkable properties. The basis $\mathcal{B}_{can}$ is however 
difficult to compute for arbitrary Dynkin diagrams.

 Let $I$ denote the set of vertices of $D$. 
Kashiwara introduced crystal operators $\tilde{e}_i, \; \tilde{f}_i, \, i \in I$ on $U_q(\nbold^+)$.
These allow to construct the crystal graph $B(\infty)$ which serves as a combinatorial skeleton 
of $\mathcal{B}_{can}$. Its vertices are the elements of $B$, and its edges are induced
by the action of crystal operators on $\mathcal{B}_{can}$. 
The crystal limit $b \mapsto b \mod q^{-1} L$
establishes a one-to-one correspondence between $\mathcal{B}_{can}$ and vertices
of $B(\infty)$, which allows to extract important combinatorial information 
from $\mathcal{B}_{can}$ to the level of $B(\infty)$. The crystal graph $B(\infty)$
may be defined by purely combinatorial means, and
provides important data for the study of finite-dimensional representations of $U_q(\gbold)$.  

The crystal limit also establishes a  one-to-one correspondence between members 
of a basis $\mathcal{P}_{\mathbf{w}_0}$ and the vertices of $B(\infty)$. 
A PBW monomial is defined by an $N$-tuple of positive integers, where $N$ is
the number of positive roots in the root system associated to $D$.
 One thus gets an indexation of $B$ by $\Nbold^N$. This is the Lusztig's parametrization
with respect to $\mathbf{w}_0$. 

Its advantage lies in the simple indexing set for $B$.
The action of a crystal operator $\tilde{e}_i$ is easy to describe 
when the reduced expression $\mathbf{w}_0$ starts with the 
simple reflection $s_i$. 
It is however difficult to give, for a fixed $\mathbf{w}_0$,  the action of all 
the operators $\tilde{e}_i, \, i \in I$ in the same time, 
due to the complexity of the passage formulas \cite{BerensteinZelevinskyInventiones} 
between Lusztig parametrizations.
This was done by Reineke \cite{Reineke}
for  reduced expressions $\mathbf{w}_0$ adapted to 
quivers $Q$ of $ADE$ type verifying a particular homological
condition  $(L)$ (detailed in section 2). 
The Hall algebra construction \cite{RingelHall} of $U_q(\nbold^+)$
allows Reineke to study the crystal operators
$\tilde{e}_i, \, i \in I$ using  the representation theory  
of finite dimensional algebras. 

Kashiwara showed \cite{KashiwaraDemazure}
that given $\mathbf{w}_0$, 
there is an elementary construction of $B(\infty)$ depending
on $\mathbf{w}_0$,
known as Kashiwara's embedding.
The vertices $B$ of $B(\infty)$ 
are indexed a set $\mathcal{S}_{\mathbf{w}_0}$ 
of specific $N$-tuples of integers, known
as string parameters. 
The action of crystal operators on $B$ is easy to describe, 
as it depends only on the Cartan matrix of $D$. 
However it is a complex problem
to describe the parameter set $\mathcal{S}_{\mathbf{w}_0}$.
 It  is the set of integer points of a polyhedral
cone $\mathcal{C}_{\mathbf{w}_0}$ \cite{LittelmannCone},
\cite{BerensteinZelevinskyInventiones}. 
A system of inequalities defining $\mathcal{C}_{\mathbf{w}_0}$
was given by 
Littelmann \cite{LittelmannCone} for particular
reduced expressions with a good structure. 
Such a set of inequalities, for arbitrary $\mathbf{w}_0$,
was given by Gleizer and Postnikov 
\cite{GleizerPostnikov} 
in $A_n$ case, and 
Berenstein and Zelevinsky 
\cite{BerensteinZelevinskyInventiones}
for all finite Dynkin types.

Any cone inequality may
be seen as $\mathbf{a} \cdot \mathbf{x} \geq 0$
where $\mathbf{a} \in \Rbold^N$ is a vector
orthogonal to the hyperplane of the inequality. 
Thus a polyhedral cone may be seen as being
defined by a finite set of vectors.  The methods of 
\cite{GleizerPostnikov} and \cite{BerensteinZelevinskyInventiones} 
construct respectively sets of vectors
$K_{\mathbf{w}_0}^{GP}, \; K_{\mathbf{w}_0}^{BZ}$ with integer coordinates, 
defining $\mathcal{C}_{\mathbf{w}_0}$. 
A move associated to a crystal operator $\tilde{e}_i$
in a given parametrization, is a vector $\vbold$
appearing as the difference between $\tbold$
and $\tilde{e}_i  \tbold$
for some $N$-tuple $\tbold$. 
We shall denote by 
$L_{\mathbf{w}_0}$ the set of all possible moves, for
all $\tilde{e}_i, \; i \in I$,
in the Lusztig parametrization with respect to $\mathbf{w}_0$.
Reineke's construction allows,
for the reduced expressions $\mathbf{w}_0$
for which it valid, 
to describe the set $L_{\mathbf{w}_0}$
in terms of the Auslander-Reiten quiver $\Gamma_Q$ of $Q$.

The main result of this article (Theorem 2.4) is that
 for $\mathbf{w}_0$ adapted to a quiver $Q$ of type $A_n$,
one has $K^{GP}_{\mathbf{w}_0}=L_{\mathbf{w}_0}$.
Thus the problem of constructing the cone 
$\mathcal{C}_{\mathbf{w}_0}$ seems to be the
the same as the one of describing action of operators
$\tilde{e}_i$ in Lusztig's parametrization for $\mathbf{w}_0$.
We conjecture the set of Lusztig moves $L_{\mathbf{w}_0}$ 
defines $\mathcal{C}_{\mathbf{w}_0}$
for reduced expression adapted to quivers $Q$ of $ADE$ type, under 
the assumption that condition $(L)$ required by Reineke on $Q$ is verified. 
We give in the last section, a $D_n$ example.

A by-product of the main theorem
is that Auslander-Reiten quivers  
allow to compute a set of defining
inequalities of $\mathcal{C}_{\mathbf{w}_0}$.
The combinatorics involved is that
of hammocks, introduced by Brenner \cite{Brenner}.
This provides an alternative to
methods given in \cite{GleizerPostnikov}, 
\cite{BerensteinZelevinskyInventiones}.
  
We would like to thank Bernard Leclerc for many helpful remarks.

\section*{2. Combinatorial models of the canonical basis}

\hspace*{5mm} Let $\Phi$ be the root system corresponding to the Dynkin diagram $D$
and $( \, , \,)$ the Cartan scalar product over
$\Rbold \Phi$.
We shall denote by $\alpha_i, i \, \in I$ the set of simple roots,
$\omega_i, i \in I$ the set of fundamental weights, 
and by $s_i, i\, \in I$ the simple reflections inside $W$.
We shall fix all through this section, a reduced expression 
$\mathbf{w}_0=s_{i_1} s_{i_2} \ldots  s_{i_N}$
of the longest element $w_0$ of $W$. This expression  induces 
the reflection ordering 
$\preccurlyeq_{\mathbf{w}_0}$ on the set of positive 
roots $\Phi^+$,  a total ordering given by 
$\beta_1=\alpha_{i_1}, \, 
\beta_2=s_{i_1} (\alpha_{i_2}), \ldots, 
\beta_N=s_{i_1} s_{i_2} \ldots s_{i_{N-1}}(\alpha_{i_N})$.

We shall use here the conventions of \cite{BerensteinZelevinskyInventiones},
concerning  the quantized enveloping algebra $U_q(\gbold)$
associated to $D$. It is generated by the set  
$e_i, f_i, k_i^{\pm1}, i \in I$
subject to relations derived from the Cartan matrix $C$.
The reader may find details of the defining relations
in \cite{BerensteinZelevinskyInventiones}, section 3.1.
We recall in particular that  $\lbrack n \rbrack_q !$ denotes the $q$-factorial of $n$,
and that the $n$-th divided power of an element $x \in U_q(\gbold)$
 is  given by $x^{(n)}=x^n / \lbrack n \rbrack_q !$.

The positive part $U_q(\nbold^+)$ is the subalgebra
generated by the $e_i, \, i  \in I$. It admits a grading by 
$Q^+=\bigoplus \Nbold \alpha_i$
obtained by putting $\deg(e_i)=\alpha_i$.
Given an arbitrary $\gamma \in Q^+$, 
the weight space $U_q(\nbold^+)_{\gamma}$ is the 
$\Cbold(q)$-vector space of elements of a degree $\gamma$.
All weight spaces of $U_q(\nbold^+)$ are finite dimensional.

The braid group acts 
on $U_q(\gbold)$ by automorphisms $T_i, \, i \in I$ (noted $T'_{i,-1}$ in part VI of
\cite{LusztigBook}).
We refer again to  \cite{BerensteinZelevinskyInventiones}
section 3 for a detailed definition.
For every $k=1, 2, \, \ldots \, , N$, 
$E_{\beta_k}:=T_{i_1} T_{i_2} \ldots T_{i_{k-1}} (E_{i_k})$
is an element of $U_q(\nbold^+)$ of weight $\beta_k$. 
Any given $N$-tuple $\mathbf{t}=(t_1, \ldots, \, t_N)$
of positive integers
defines the \textbf{PBW monomial}
$$ p_{\mathbf{w}_0}(\mathbf{t}):=
E_{\beta_1}^{(t_1)} E_{\beta_2}^{(t_2)} 
\ldots E_{\beta_N}^{(t_N)}.$$
The set of all such monomials, 
$\mathcal{P}_{\mathbf{w}_0}:=
\{ p_{\mathbf{w}_0}(\mathbf{t}) \mid \mathbf{t} \in \Nbold^N \}$, 
forms the \textbf{PBW-basis} of $U_q(\nbold^+)$ associated 
to the reduced expression $\mathbf{w}_0$.

\textbf{Theorem 2.1} \cite{LusztigBase}

\emph{For every monomial $p_{\mathbf{w}_0}(\mathbf{t})$, there is one and only one $b \in \mathcal{B}_{can}$
such that \linebreak[4] $b=p_{\mathbf{w}_0}(\mathbf{t}) \, \mod \, q^{-1} \mathcal{L}$.} 

The crystal operators $\tilde{e}_i, \tilde{f}_i, \, i \in I$ where introduced
by Kashiwara \cite{KashiwaraBase} for the \linebreak[4] negative part
$U_q(\nbold^-)$ of $U_q(\gbold)$. As $U_q(\nbold^-)$ and
$U_q(\nbold^+)$ are isomorphic as algebras, this construction
may be carried over to $U_q(\nbold^+)$.

Given $i \in I$, there is a locally nilpotent action $\theta_i$  over $U_q(\nbold^+)$ 
defined by  :
$$ 
\theta_i(1)=0, \; \;  \forall x \in U_q(\nbold^+) \; : \; \theta_i(e_j x)=q^{(\alpha_i, \alpha_j)} e_j \theta_i(x)+\delta_{i,j} x.$$
One has $U_q(\nbold^+)=\underset{n \in \Nbold}{\bigoplus} e_i^{(n)} \ker \theta_i$,
and $\ker \theta_i$ is compatible with the weight graduation of
$U_q(\nbold^+)$. One chooses a weight vector basis $\Xi_i$ of $\ker \theta_i$ , 
and defines for each $\vbold \in \Xi_i$ and
any $n \in \Nbold$ :
$$\begin{array}{rl}
\tilde{e}_i (e_i^{(n)} v) := & e_i^{(n+1)} v, \\
\tilde{f}_i (e_i^{(n)} v) : = & \left\{
\begin{array}{ll} e_i^{(n-1)} v & \mbox{if} \; n \geq 1, \\
                        0                 & \mbox{if} \; n=0.                                    
\end{array}
\right.
 \end{array}   
$$             
This leads to well defined operators $\tilde{e}_i, \, \tilde{f}_i$
over $U_q(\nbold^+)$, which do not depend on the initial choice of $\Xi_i$.

The crystal operators $\tilde{e}_i, \tilde{f}_i, \, i \in I$ 
preserve $\mathcal{L}$, and hence induce an action
over $\mathcal{L}/q^{-1} \mathcal{L}$.
A key feature of the canonical basis is its good behaviour under 
this action.
For any $b \in \mathcal{B}_{can}$,
$\tilde{e}_i b = b' \mod q^{-1} \mathcal{L}$,
and $\tilde{f}_i b$ is either $ 0 \mod q^{-1} \mathcal{L}$ or
\linebreak[4]
$\tilde{f}_i b=b'' \mod q^{-1} \mathcal{L}$, where
$b'$, $b''$ are other elements of $\mathcal{B}_{can}$.
We see the image $B$ of $\mathcal{B}_{can}$ inside
$\mathcal{L}/q^{-1} \mathcal{L}$ becomes endowed
with a structure of a coloured graph, the \linebreak[4] arrows
being valuated by the operators $\tilde{e}_i, \tilde{f}_i, \, i \in I$.
This is the crystal graph $B(\infty)$ of Kashiwara \cite{KashiwaraBase}.

Our discussion in the introduction and Theorem 2.1 above lead
to a one-to-one correspondence 
$\varphi_{\mathbf{w}_0} \, : \, \Nbold^N \longrightarrow B, \; \; 
\tbold \mapsto \pbold_{\mathbf{w}_0}(\tbold) \mod q^{-1} \mathcal{L}$.  
This is \textbf{Lusztig's parametrization} with respect to $\mathbf{w}_0$.
Under this identification, we may consider the crystal operators  $\tilde{e}_i, 
\tilde{f}_i$ as acting on $\Nbold^N$. 
We shall call a vector $\lbold \in \Zbold^N$ 
a \textbf{Lusztig move of type} $i$ with respect 
to $\mathbf{w}_0$, if there exists $\tbold \in \Nbold^N$
such that $\tilde{e}_i \tbold=\tbold+\lbold$.
Recall $L_{\mathbf{w}_0}$ denotes the set
of all possible Lusztig moves, for all types $i \in I$.

\textbf{Example : } 

Consider $A_2$ case, $\wbold_0=s_1 s_2 s_1$.
  For a given $\tbold=(t_1, \, t_2, \, t_3)$ one has :
$$ \begin{array}{rl}
\tilde{e}_1 (t_1, \, t_2, \, t_3)& =(t_1+1,\, t_2, \, t_3) \\
\tilde{e}_2 (t_1, \, t_2, \, t_3)& =
\left\{ 
\begin{array}{ll}
(t_1-1,\, t_2+1, \, t_3) & \mbox{if} \; t_1 > t_3 \\
(t_1, \, t_2, \, t_3+1) & \mbox{if} \; t_1 \leq t_3.
\end{array} \right.
\end{array}$$

One sees there is only one Lusztig move of type 1,
$\mathbf{l}_1=(1,0,0)$ and two Lusztig moves
of type 2, $\mathbf{l}_2=(-1,1,0), \;
\mathbf{l}_3=(0,0,1).$ We get
$L_{\mathbf{w}_0}=\{ (1,0,0), \, (-1,1,0), \, (0,0,1) \}$.
\medskip

Let $Q$ be a fixed quiver obtained 
by orienting the Dynkin diagram $D$.
Following $\cite{BernsteinGelfandPonomarev}$, we call
a vertex $i$ of Q a sink, if there are only 
arrows entering it. We denote in that case
by $s_i Q$ the quiver obtained by reversing
the arrows whose end is $i$, into arrows exiting
$i$, thus transforming $i$ into a source. 
A reduced expression
$\mathbf{w}_0=s_{i_1} s_{i_2} \ldots s_{i_N}$
is adapted to Q if and only if $i_1$ is a sink of $Q$,
$i_2$ a sink of $s_{i_1} Q$, $i_3$ a sink of 
$s_{i_2} s_{i_1} Q$ and so on. Such an expression always
exists for a given $Q$.

Let us  denote  by $\Cbold Q$ the path algebra of $Q$ over $\Cbold$.
The  category $\mbox{mod}\; \Cbold Q$ of 
finite dimensional left modules has simple objects
$S_i$ which are indexed by $I$. 
We shall say, following Reineke \cite{Reineke}, that
 the quiver $Q$ verifies \textbf{condition} $(L)$, if for every 
indecomposable module
$X \in \mbox{mod} \; \Cbold Q$,  and every 
$i \in I$, one has $\dim \mbox{Hom}(X,S_i) \leq 1. $
This condition is verified for any quiver of type $A_n$, 
and at least one quiver of each of the types $D_n, \, E_6, \, E_7$
(\cite{Reineke} Appendix).
Under this condition,
the action of crystal operators in a Lusztig parametrization
with respect to  $\mathbf{w}_0$  adapted to $Q$ 
may be described in terms of the category $\mbox{mod}\; \Cbold Q$,
and the set $L_{\mathbf{w}_0}$ may be obtained \cite{Reineke}. 
We postpone the details to section 3.

We continue to fix the same $\wbold_0 =s_{i_1}, \ldots s_{i_N}$.
We refer the reader to 
Kashiwara \cite{KashiwaraDemazure}, and \cite{JosephBook} Chapter 5
for details on crystal theory. 
Kashiwara's elementary construction of $B(\infty)$
uses the Cartan matrix $C=(c_{i,j})$
in order to define operators $\tilde{e}_i, \tilde{f}_i, \, i \in I$
acting over $\Nbold^N$ as below (\cite{JosephBook} 5.2.5, 6.1.15) : 

Fix $\abold=(a_1, \, \ldots, \, a_N) \in \Nbold^N$.  For 
$k=1, 2, \ldots N$, define  $r_k:=a_k+\underset{1 \leq j <k}{\sum} c_{i_j,i_k} a_j$.
Given $i \in I$, consider $\xi_i=\max_{i_k=i} r_k$. Let $k_1$
be the first position where this maximum is attained, $k_2$ the last. Then :
$$ \begin{array}{rl}
\tilde{e}_i (\abold)& = (a_1, \, \ldots, \, a_{k_2-1}, \, a_{k_2}+1, \, a_{k_2+1}, \, \ldots, \, a_N). \\
\tilde{f}_i (\abold)& =  
\left\{
\begin{array}{cl} 
(a_1, \, \ldots, \, a_{k_1-1}, \, a_{k_1}-1, \, a_{k_1+1}, \, \ldots, \, a_N) & \mbox{if} \, a_{k_1} \geq 1, \\
0 & \mbox{if} \, a_{k_1}=0.
\end{array}
\right. 
\end{array}$$

\textbf{Theorem 2.2} : \emph{Kashiwara's embedding} \cite{KashiwaraBase}

\emph{
Let $\mathcal{N}$ denote $\Nbold^N$ with the action
of operators $\tilde{e}_i, \, \tilde{f}_i, \, i \in I$ given above.
\begin{itemize}
\item[a)] There is an embedding $\psi_{\mathbf{w}_0} \; B
\hookrightarrow \Nbold^N$, sending the graph $B(\infty)$ isomorphicaly onto 
the subgraph of $\mathcal{N}$ generated 
out of  $\ubold_{\infty}:=(0,0, \ldots 0)$ by applying the operators $\tilde{e}_i, i \in I$.
\item[b)] The image of $B(\infty)$ consitsts of those 
$\abold=(a_1, a_2 , \ldots , a_N) \in \Nbold^N$ verifiying
$$  \forall k=1 \ldots N \; :  \hspace*{5mm}
\tilde{f}_{i_k} (\tilde{e}_{i_{k-1}}^{a_{k-1}} 
\tilde{e}_{i_{k-2}}^{a_{k-2}} \ldots \tilde{e}_{i_1}^{a_1} \ubold_{\infty})=0.$$
\end{itemize}}
  
Elements $\mathbf{a}:=(a_1, \, a_2, \ldots a_N)$ in b)
above are called \textbf{string parameters} \cite{BerensteinZelevinskyStrings}. 
The parametrization  the set $B$ of vertices $B(\infty)$ obtained
through Theorem 2.2  is \linebreak[4] \textbf{Kashiwara's parametrization} with 
respect to $\textbf{w}_0$. 
Recall we  denote its indexing set 
$\mbox{Im} \, \psi_{\mathbf{w}_0}$ by $\mathcal{S}_{\mathbf{w}_0}$.

\textbf{Remarks}

i) Kashiwara works with the negative part $U_q(\nbold^-)$
of $U_q(\gbold)$. The definition above is the transfer of his construction
to $U_q(\nbold^+)$,  which amounts to exchanging the roles 
of $\tilde{e}_i$ and $\tilde{f}_i$ at the level of $B(\infty)$.

ii) The definition of Kashiwara's embedding imposes a reversal
of order in the \linebreak[4] definition of a string, namely
$(a_1, a_2, \ldots a_N)$ in our convention, corresponds
to \linebreak[4]
$(a_N, a_{N-1}, \ldots a_1)$ in \cite{BerensteinZelevinskyStrings}.

\textbf{Theorem 2.3} \cite{LittelmannCone}, 
\cite{GleizerPostnikov}, 
\cite{BerensteinZelevinskyInventiones}  

\emph{
The set $\mathcal{S}_{\mathbf{w}_0}$ is the set of integer 
points of a polyhedral cone $\mathcal{C}_{\mathbf{w}_0}$, that is, there exists 
a finite set of vectors 
$K_{\mathbf{w}_0} \subset \Zbold^N$ such that 
$$ \mathcal{S}_{\mathbf{w}_0}=\{ \mathbf{a} \in \Nbold^N \mid
\forall \mathbf{k} \in K_{\mathbf{w}_0},  \;  
\mathbf{a} \cdot \mathbf{k} \geq 0 \}.$$}

\textbf{Example}

Consider type $A_2$, and $\mathbf{w}_0=s_1 s_2 s_1$.
It is easy to compute the image of Kashiwara's embedding  
directly out of its definition above.
One obtains the well known result 
$$ \mathcal{C}_{\mathbf{w}_0}=
\{ (a_1,\, a_2 , \, a_3) \mid 
0 \leq a_1 \leq a_2, 0 \leq a_3 \}.$$
One may choose as a defining set for $\mathcal{C}_{\mathbf{w}_0}$,  
the set $K_{\mathbf{w}_0}=\{ (1,0,0), \, (-1,1,0), \, (0,0,1) \}$, which
is equal to both $K_{\mathbf{w}_0}^{GP}$ and $K_{\mathbf{w}_0}^{BZ}$.
\medskip

\textbf{Main Theorem 2.4}

\emph{
Let  $\mathbf{w}_0$ be a reduced expression adapted to
a quiver of type $A_n$, and $K^{GP}_{\mathbf{w}_0}$
the set given by Gleizer and Postnikov (\cite{GleizerPostnikov} section 5). Then} 
$$  K^{GP}_{\mathbf{w}_0}=L_{\mathbf{w}_0}.$$

Let us observe that in Lusztig's parametrization,
the parameter set is the set of integer points of the 
cone $(\Rbold^+)^N$, consisting of vectors with positive coordinates. 
This cone may be seen as being defined by the natural basis 
$E=\{ \ebold_1, \ebold_2 \ldots \ebold_N \}$ of $\Rbold^N$.
It is easy to see, by the definition of Kashiwara's embedding, 
that the set $E$ is the set of
of moves of crystal operators $\tilde{e}_i, \, i \in I$ 
in Kashiwara's parametrization according to $\mathbf{w}_0$.
We have therefore a full symmetry between Lusztig's and Kashiwara's parametrizations, 
the set of vectors defining the parameters set in one picture being 
equal to the set of crystal operators 
moves in the other.

\textbf{Conjecture} 

\emph{
Let $Q$ be a quiver of type $ADE$ satisfying Reineke's condition
$(L)$. Let $\mathbf{w}_0$ be adapted to it. Then 
$L_{\mathbf{w}_0}$ is a 
defining set for $\mathcal{C}_{\mathbf{w}_0}$.}

We state this conjecture on the basis of
some computer testing, using the sets $K^{BZ}_{\mathbf{w}_0}$.
We give in section 7 a detailed $D_4$ example.
The conjecture might be valid in a larger
scope, even beyond reduced expressions adapted to quivers.
However one faces a breakdown of 
many nice properties, enjoyed by reduced expressions 
adapted to quivers verifying condition $(L)$.

\section*{3. Auslander-Reiten quivers and Lusztig's moves}

\hspace*{5mm} Fix $Q$ a quiver of type $ADE$ satisfying condition
$(L)$, $\mathbf{w}_0$ adapted to it and 
$\beta_1, \, \beta_2, \, \ldots, \, \beta_N$ the reflection 
ordering it defines.
 The  category $\mbox{mod}\; \Cbold Q$ is \linebreak[4] equivalent
to that of finite dimensional representations of $Q$.
A module $M$ in $\mod \; \Cbold Q$ may be seen
as a family $(V_i)_{i \in I}$ of finite dimensional $\Cbold$-vector
spaces, together with \linebreak[4] linear mappings $f_{i,j} \; : \; V_i \longrightarrow V_j$
corresponding to the arrows $i \longrightarrow j$ of $Q$.
The dimension vector of $M$ is the element of $Q^+$ 
given by $\dbold_M:=\underset{i=1}{\overset{n}{\sum}} (\dim V_i) \alpha_i$.
A simple object  $S_i$ has a  dimension vector equal to $\alpha_i$.

Let $\mbox{Ind} \, Q$ denote the set of isomorphism
classes of indecomposable objects of $\mod \; \Cbold Q$.
The theorem of Gabriel states that for each $\beta \in \Phi^+$, 
there is a unique class $\lbrack M \rbrack \in \mbox{Ind} \, Q$
with $\dbold_M=\beta$, and that all indecomposable objects of
$\mod \, \Cbold Q$ are obtained this way. There is therefore a one-to-one
correspondence between $\mbox{Ind} \, Q$ and $\Phi^+$, 
and we shall  denote by $\lbrack \beta \rbrack$
 the  class $\lbrack M \rbrack$ in $\mbox{Ind} \, Q$ whose dimension vector is $\beta$.

The \textbf{Auslander-Reiten quiver} 
$\Gamma_{Q}$ (\cite{AuslanderReitenSmalo} Chapter VII, \cite{GabrielARquivers} section 6)
has as set of vertices $\mbox{Ind} \, Q$, 
and its arrows are irreducible morphisms between objects of $\mbox{Ind} \, Q$.
It has a rigid mesh structure, as given in \cite{GabrielARquivers} Figure 13, page 49.
The quiver $\Gamma_Q$ is endowed with the translation 
$\tau$ (\cite{AuslanderReitenSmalo} page 225) which
sends non-projective modules of $\mbox{Ind} \,  Q$ 
onto non-injective modules of $\mbox{Ind} \, Q$.
 The translation $\tau$ stratifies $\Gamma_Q$ into levels. 
The $i^{th}$ level 
is  the orbit under $\tau$ of the injective envelope of 
$S_{i^*}$, where $*$ denotes the Dynkin diagram automorphism
induced by $w_0$. This level ends  in the projective cover
of $S_i$.

There is a natural order on the vertices of $\Gamma_Q$, given by
$\lbrack \beta_1 \rbrack \leq_Q \lbrack \beta_2 \rbrack$ if
and only if there 
is a path from  $\lbrack \beta_1 \rbrack$ to 
$\lbrack \beta_2 \rbrack$ 
in $\Gamma_Q$.
This induces a partial order
$\preccurlyeq_Q$ on $\Phi^+$ by putting
$\beta_1 \preccurlyeq_Q \beta_2$ 
whenever 
$\lbrack \beta_1 \rbrack \leq_Q \lbrack \beta_2 \rbrack$.
 The reflection ordering ordering 
$\preccurlyeq_{\mathbf{w}_0}$
is then a linear refinement of $\preccurlyeq_Q$.

The path algebra $\Cbold Q$ is an hereditary algebra. The Euler-Poincar\'{e} characteristic  
$ \langle M_1 , \, M_2 \rangle:=\dim \mbox{Hom}(M_1,M_2)
-\dim \mbox{Ext}^1 (M_1,M_2)$ depends only on the
dimension vectors $\dbold_{M_1}, \dbold_{M_2}$ 
 of $M_1$ and $M_2$. One has
$\langle \lbrack \beta_1 \rbrack, \lbrack \beta_2 \rbrack \rangle=
(\beta_1, \, \beta_2)_R$, where $(, )_R$
is the \textbf{Ringel form} upon 
the Euclidean space $\Rbold \Phi$. 
The matrix $R=(r_{i,j})$ of this form,  
in the basis of the simple roots, is given by
$$ r_{i,j}:=(\alpha_i, \, \alpha_j)_R=
\left\{ \begin{array}{cl} 
1 & \mbox{if  $i=j$,} \\
-1 & \mbox{if  $i \longrightarrow j$ is in Q,} \\
0 & \mbox{otherwise.} 
\end{array} \right. $$

\textbf{Theorem 3.1} \cite{RingelPBW}
\emph{
\begin{itemize}
\item[i)] Suppose $\beta_1 \preccurlyeq_Q \beta_2$ then
$\dim Ext^1 (\lbrack \beta_1 \rbrack, 
\lbrack \beta_2 \rbrack)=0$ 
and therefore \hfill \break $\dim Hom(\lbrack \beta_1 \rbrack, \,
\lbrack \beta_2 \rbrack) = (\beta_1,\, \beta_2)_R$.
\item[ii)] Suppose $\beta_1 \succcurlyeq_Q \beta_2$ then
$\dim Hom(\lbrack \beta_1 \rbrack, \, 
\lbrack \beta_2 \rbrack)=0$ and therefore \hfill \break
$\dim Ext^1 (\lbrack \beta_1 \rbrack, 
\lbrack \beta_2 \rbrack)=-(\beta_1,\, \beta_2)_R$.
\end{itemize}}

Theorem 3.1 reduces the testing of condition $(L)$ 
for a quiver $Q$, to computations in terms
of $\Phi^+$.

We refer to Section 2 of \cite{Reineke} for a concise description
of the Hall algebra construction of $U_q(\nbold^+)$
and its link to PBW bases. 
The $\Cbold(q)$ vector space with formal base vectors $\ubold_{\lbrack M \rbrack}$
indexed by isomorphism classes of $\mod \, \Cbold Q$,
may be endowed with a product linked to the module structure. 
This defines the Hall algebra $\mathcal{H}(Q)$. Ringel's main 
theorem \cite{RingelHall} states that sending the generators 
$e_i$ to $\ubold_{\lbrack S_i \rbrack}$,
establishes an isomorphism 
$\eta_Q \; : \; U_q(\nbold^+) \stackrel{\sim}{\longrightarrow} \mathcal{H}(Q)$.

Let $ \lbrack M \rbrack=\underset{j=1}{\overset{N}{\bigoplus}}  
\lbrack \beta_j \rbrack^{\oplus t_j}$ be an isoclass 
with multiplicities of indecomposables
given by $\tbold_M:=(t_1, t_2, \, \ldots \, t_N)$. 
The PBW  monomial $\pbold_{\mathbf{w}_0}(\tbold_M)$, with 
$\mathbf{w}_0$ adapted to $Q$ is recovered,  up
to a multiplication by a well defined power of $q$, 
as the inverse image under $\eta_Q$ of $\ubold_{\lbrack M \rbrack}$. 
Crystal operators in the Lusztig parametrization of 
$\mathbf{w}_0$ may therefore be seen
as acting upon isomorphism classes of 
$\mbox{mod} \; \Cbold Q$. 
One has $\tilde{e}_i \lbrack M_1 \rbrack =\lbrack M_2 \rbrack$
if and only if \linebreak[4] 
$\tilde{e}_i \pbold_{\mathbf{w}_0}(\tbold_{M_1})=
\pbold_{\mathbf{w}_0}(\tbold_{M_2 }) 
\mod \, q^{-1} \mathcal{L}$.

Let us fix $i \in I$. The description of the action
of $\tilde{e}_i$ is given in terms of the set 
(\cite{Reineke} page 711) :
$$ P_i(Q):=\{ \lbrack X \rbrack \in \mbox{Ind} Q \mid 
\dim  Hom (X , S_i) > 0 \}.$$

The set $P_i(Q)$ has a poset structure $\lbrack X \rbrack \leq \lbrack Y \rbrack$ 
whenever there is a path from $\lbrack X \rbrack$ to  $\lbrack Y \rbrack$
inside $P_i(Q)$. It is the same as the order induced by $\leq_Q$
(\cite{Reineke} Proposition 4.3). Recall an \textbf{antichain} 
$A$ of a $P_i(Q)$ is a set of mutually non-comparable elements.
It defines the \textbf{order ideal}
$J(A):=\{ \lbrack X \rbrack \in P_i(Q) \mid 
\exists \lbrack Y \rbrack \in A, \lbrack X \rbrack \leq \lbrack Y \rbrack \}$.
The correspondence $A \mapsto J(A)$ is one-to-one, 
and inclusion between order ideals induces
a natural poset structure upon the set $\mathcal{A}_i(Q)$
of all antichains of $P_i(Q)$. 

Given  $A \in \mathcal{A}_i(Q)$, let $C_A$ 
be the set of minimal elements of 
$P_i(Q) \backslash J(A)$, and define 
$$ \begin{array}{rl}
\lbrack V_{A} \rbrack:= & 
\underset{\lbrack M \rbrack \in A}{\oplus} \lbrack M \rbrack; \\
\lbrack U_{A} \rbrack := &
\underset{\lbrack M \rbrack \in C_A}{\oplus} \lbrack \tau M \rbrack.
\end{array} $$   
 
The Lusztig move corresponding to $A$ is then
$\mathbf{l}_A:=\tbold_{V_A}-\tbold_{U_A}$. 

Each $A \in \mathcal{A}_i(Q)$ also defines 
a function $F_A \; : \; \Nbold^{N} \longrightarrow \Zbold$
given by 
$$ F_{A}(\tbold):=
\sum_{X \in J(A)}
\Delta_X(\tbold)$$
where for $X=\lbrack \beta_k \rbrack$ with
$\tau X=\lbrack \beta_{k'} \rbrack$, 
$\Delta_X(\tbold)=t_k-t_{k'}$ (with the convention
that the second term is $0$ if translation
is not defined on $X$).

\textbf{Theorem 3.2}\emph{ (\cite{Reineke} Theorem 7.1)}

\emph{
Let $\lbrack M \rbrack$ be an isomorphism
class  of $\mbox{mod} \, \Cbold Q$. Put
$\zeta_i(\lbrack M \rbrack):=\max_{A \in \mathcal{A}_i(Q)} 
F_{A}(\tbold_M)$. Then the subset 
$\{ A \mid F_{A}(\tbold_M)=\zeta_i(\lbrack M \rbrack) \}$
of $\mathcal{A}_i(Q)$ 
admits a unique
maximal element $A_{max}$.
There is an isoclass $\lbrack X \rbrack$ such that
 $\lbrack M \rbrack  =\lbrack X \rbrack 
  \oplus  \lbrack U_{A_{max}} \rbrack$, and the action
of $\tilde{e}_i$ is given then by
  $\tilde{e}_i \lbrack M \rbrack =\lbrack X \rbrack
   \oplus \lbrack V_{A_{max}} \rbrack$.}

In terms of Lusztig parameters, 
$\tilde{e}_i \tbold_M=\tbold_M+\lbold_{A_{max}}$.

\textbf{Corollary 3.3} 

\emph{
Consider the Lusztig parametrization
corresponding to $\mathbf{w}_0$. 
The set  of Lusztig moves of type $i$ is
given by
$$L_{\mathbf{w}_0}^{(i)}=
\{ \mathbf{l}_A \mid A \in \mathcal{A}_i(Q) \}.$$}  

\textbf{Proof} 

In view of Theorem 3.2,  
given $A \in \mathcal{A}_i(Q)$, 
we need to exhibit a module  
on which $\tilde{e}_i$ acts according
to the move defined by $A$. Such a  module 
is $U_A$.
One has $\Delta_X(\tbold_{U_A}) \geq 0$
for any $\lbrack X \rbrack \in J(A)$,
where as  $\Delta_X(\tbold_{U_A}) = -1$
for any $\lbrack X \rbrack \in C_A$. These
two properties ensure that for $U_A$, $A_{max}$
of Theorem 3.2 is $A \; .\Box$
\medskip

\textbf{Example}

Let $Q$ be the quiver
$\stackrel{1}{\cdot} \longleftarrow
\stackrel{2}{\cdot} \longrightarrow \stackrel{3}{\cdot}$
of type $A_3$, with adapted reduced expression \linebreak[4]
$\mathbf{w}_0=s_1 s_3 s_2 s_1 s_3 s_2$. 
The Auslander-Reiten quiver $\Gamma_Q$ is 
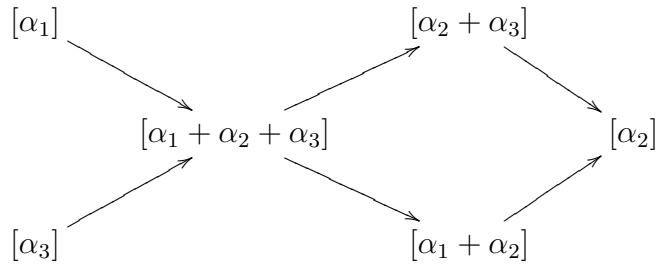
\[
\xymatrix{
\lbrack \alpha_1 \rbrack \ar[rd] & & 
\lbrack \alpha_2+\alpha_3 \rbrack \ar[rd] & \\
& \lbrack \alpha_1+\alpha_2+\alpha_3 \rbrack \ar[ru] \ar[rd] & 
& \lbrack \alpha_2 \rbrack \\
\lbrack \alpha_3 \rbrack \ar[ru] & & 
\lbrack \alpha_1+\alpha_2 \rbrack \ar[ru] & } 
\] 

Let us study
the action of $\tilde{e}_2$. Using $( \, , \,)_R$,
and Theorem 3.1, one gets \linebreak[4]
$P_2(Q)=\{ \lbrack \alpha_2 \rbrack, 
\lbrack \alpha_1+\alpha_2 \rbrack,
\lbrack \alpha_2+\alpha_3 \rbrack, 
\lbrack \alpha_1+\alpha_2+\alpha_3 \rbrack \}$. 
$\mathcal{A}_2(Q)$ consists in 5 antichains : \medskip
\hfill \break
\begin{tabular}{|l|c|c|c|}
Antichain $A$ & $U_A$ & $V_A$ & $\lbold_A$ \\
\hline
$A_1=\{ \lbrack \alpha_1+\alpha_2+\alpha_3 \rbrack \}$ &
$\lbrack \alpha_1 \rbrack \oplus \lbrack \alpha_3 \rbrack$ &
$\lbrack \alpha_1+\alpha_2+\alpha_3 \rbrack$ &
$(-1,-1,1,0,0,0)$ \\
$A_2=\{ \lbrack \alpha_2+\alpha_3 \rbrack \}$ &
$\lbrack \alpha_3 \rbrack$ &
$\lbrack \alpha_2+\alpha_3 \rbrack$ &
$(0,-1,0,1,0,0)$ \\
$A_3=\{ \lbrack \alpha_1+\alpha_2 \rbrack \}$ &
$\lbrack \alpha_1 \rbrack$ &
$\lbrack \alpha_1+\alpha_2 \rbrack$ &
$(-1,0,0,0,1,0)$ \\
$A_4=\{ \lbrack \alpha_1+\alpha_2 \rbrack, \,  
\lbrack \alpha_2+\alpha_3 \rbrack \}$ &
$\lbrack \alpha_1+\alpha_2+\alpha_3 \rbrack$ &
$\lbrack \alpha_1+\alpha_2 \rbrack \oplus
\lbrack \alpha_2+\alpha_3 \rbrack$ &
$(0,0,-1,1,1,0)$ \\
$A_5=\{ \lbrack \alpha_2 \rbrack \}$ &
 $0$ &
$\lbrack \alpha_2 \rbrack$ &
$(0,0,0,0,0,1)$ \\
\hline
\end{tabular}
\medskip

The set $\mathcal{A}_2(Q)$ has the following poset structure
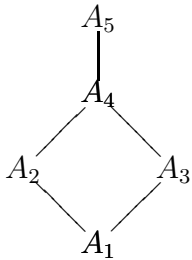
\[ \begin{xy}
(20,10)*{A_1};
(10,20)*{A_2} **\dir{-};
(20,30)*{A_4} **\dir{-};
(20,10)*{A_1};
(30,20)*{A_3} **\dir{-};
(20,30)*{A_4} **\dir{-};
(20,40)*{A_5} **\dir{-};
\end{xy}
\]

Considering the structure of $\Gamma_Q$ and 
the order ideals defined by members 
of $\mathcal{A}_2(Q)$, we get the $F$-functions
data.
\begin{center}
\begin{tabular}{c|l|l}
Antichain & $J(A)$ & $F_A$ \\
\hline
$A_1$ & $\{ \lbrack \alpha_1+\alpha_2+\alpha_3 \rbrack \}$ & $t_3$ \\
$A_2$ & $\{ \lbrack \alpha_1+\alpha_2+\alpha_3 \rbrack, \,  
\lbrack \alpha_2+\alpha_3 \rbrack  \}$ & $t_3+(t_4-t_1)$ \\
$A_3$ & $\{ \lbrack \alpha_1+\alpha_2+\alpha_3 \rbrack, \,
\lbrack \alpha_1+\alpha_2 \rbrack  \}$ & $t_3+(t_5-t_2)$ \\
$A_4$ & $\{ \lbrack \alpha_1+\alpha_2+\alpha_3 \rbrack, \,
\lbrack \alpha_1+\alpha_2 \rbrack, \, 
\lbrack \alpha_2+\alpha_3 \rbrack  \}$ & $t_3+(t_4-t_1)+(t_5-t_2)$ \\
$A_5$ & $P_2(Q)$ & $(t_4-t_1)+(t_5-t_2)+t_6$ \\
\hline
\end{tabular}
\end{center} 

Take $\tbold_M=(3,2,1,1,2,0)$. One has 
$t_3=1, \, t_4-t_2=-2, \, t_5-t_2=0, \, t_6=0.$ Replacing
each antichain $A$ by the value of $F_{A}(\tbold_M)$
gives the following diagram
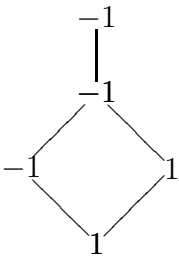
\[ \begin{xy}
(20,10)*{1};
(10,20)*{-1} **\dir{-};
(20,30)*{-1} **\dir{-};
(20,10)*{1};
(30,20)*{1} **\dir{-};
(20,30)*{-1} **\dir{-};
(20,40)*{-1} **\dir{-};
\end{xy}
\]
This pinpoints $A_3$ as $A_{max}$ for $\lbrack M \rbrack$. Thus
$\tilde{e}_2 \tbold_M=\tbold_M+\lbold_{A_3}=(2,2,1,1,3,0)$.

\section*{4. Wiring diagrams and string cones}

\hspace*{5mm} We shall restrict ourselves from now on 
to $A_n$ type. The positive roots are \linebreak[4]
$\beta=\alpha_i + \alpha_{i+1}+ \ldots +\alpha_j , \; 
1 \leq i \leq j \leq n$. They are in a one-to-one \linebreak[4] correspondence
with pairs $(i,j), \; 1 \leq i <  j \leq n+1$,
$\beta$ above being sent to $(i,j+1)$. \linebreak[4] The
fundamental representation $E(\omega_1)$ of $\slalg_{n+1}$ has weights
given by \linebreak[4] $\nu_j=-\omega_{j-1}+\omega_j, \; j=1 \ldots n+1$ 
(with the convention that $\omega_0=\omega_{n+1}=0$).
The weight $\nu_j$ may be seen as the weight of a one-box 
Young tableau $\fbox{$j$}$. The weights of 
fundamental representations $E(\omega_k)$, for $\; 2 \leq k \leq n$ correspond 
to strictly increasing column tableaux of size $k$.
The weight of a tableau is the sum of weights of its boxes.
We shall therefore identify these weights with
$k$-tuples $1 < j_1 < j_2 < \ldots < j_k \leq n+1$. 
The action of $W$ on roots or on weights 
identifies with that of the symmetric group
$\gothicS_{n+1}$ on the respective multi-indices
we considered.  

The \textbf{wiring diagram} $\mathcal{WD}(\mathbf{w}_0)$ 
of a reduced expression
$\mathbf{w}_0=s_{i_1} s_{i_2} \ldots s_{i_N}$ consists
in encoding $\mathbf{w}_0$ as an arrangement of pseudo-lines 
$L_1, \ldots , L_{n+1}$ drawn inside a vertical strip
of $\Rbold^2$. The respective
crossing points of these strands 
occur in levels according to 
the  indices of $\mathbf{w}_0$. Figure 4.1
gives a self-explaining example of the procedure in type
$A_3$.

Each strand $L_i$ crosses another strand $L_j$ once, and
only once. The order of the strands $L_i, \;  i=1, \ldots n+1$
gets inverted while following $\mathcal{WD}(\mathbf{w}_0)$ 
from left to right. Sending the crossing $v_{i,j} \; (i<j)$ 
of lines $L_i$ and $L_j$ onto the couple $(i,j)$, establishes
a one-to-one order preserving correspondence between 
these crossings enumerated from left to right and
the reflection ordering $\preccurlyeq_{\mathbf{w}_0}$
of $\Phi^+$. If $\beta_k$ is the $k^{th}$ root in that
order, with $i_k=i$, then the $k^{th}$ crossing
of $\mathcal{WD}(\mathbf{w}_0)$, that  
we shall denote by $v_{\beta_k}$,
is on level $i$.

\[
\begin{xy} 
(20,60)*{\cdot}; 
(25,60)*{} **\dir{-};
(30,55)*{\cdot} **\dir{-};
(35,50)*{} **\dir{-};
(45,50)*{} **\dir{-};
(50,45)*{\cdot} **\dir{-};
(55,40)*{} **\dir{-};
(65,40)*{} **\dir{-};
(70,35)*{\cdot} **\dir{-};
(75,30)*{} **\dir{-};
(85,30)*{} **\dir{-};
(90,30)*{\cdot} **\dir{-};
(20,50)*{\cdot};
(25,50)*{} **\dir{-};
(30,55)*{\cdot} **\dir{-};
(35,60)*{} **\dir{-};
(55,60)*{} **\dir{-};
(60,55)*{\cdot} **\dir{-};
(65,50)*{} **\dir{-};
(75,50)*{} **\dir{-};
(80,45)*{\cdot} **\dir{-};
(85,40)*{} **\dir{-};
(90,40)*{\cdot} **\dir{-};
(20,40)*{\cdot}; 
(35,40)*{} **\dir{-};
(40,35)*{\cdot} **\dir{-};
(45,30)*{} **\dir{-};
(65,30)*{} **\dir{-};
(70,35)*{\cdot} **\dir{-};
(75,40)*{} **\dir{-};
(80,45)*{\cdot} **\dir{-};
(85,50)*{} **\dir{-};
(90,50)*{\cdot} **\dir{-};
(20,30)*{\cdot};
(35,30)*{} **\dir{-};
(40,35)*{\cdot} **\dir{-};
(45,40)*{} **\dir{-};
(50,45)*{\cdot} **\dir{-};
(55,50)*{} **\dir{-};
(60,55)*{\cdot} **\dir{-};
(65,60)*{} **\dir{-};
(90,60)*{\cdot} **\dir{-};
(30,55)*{\cdot}; (40,35)*{\cdot}; 
(50,45)*{\cdot}; (60,55)*{\cdot}; 
(70,35)*{\cdot}; (80,45)*{\cdot};
(30,60)*{v_{12}}; (40,40)*{v_{34}}; (50,50)*{v_{14}};
(60,60)*{v_{24}}; (70,40)*{v_{13}}; (80,50)*{v_{23}};
(10,60)*{L_1}; (10,50)*{L_2}; (10,40)*{L_3}; (10,30)*{L_4};
(17,60)*{l_1}; (17,50)*{l_2}; (17,40)*{l_3}; (17,30)*{l_4};
(93,60)*{r_4}; (93,50)*{r_3}; (93,40)*{r_2}; (93,30)*{r_1};
(30,25)*{1}; (40,25)*{3}; (50,25)*{2}; 
(60,25)*{1}; (70,25)*{3}; (80,25)*{2}; 
(30,27)*{}; (30,53)*{} **\dir{--};
(40,27)*{}; (40,33)*{} **\dir{--};
(50,27)*{}; (50,43)*{} **\dir{--};
(60,27)*{}; (60,53)*{} **\dir{--};
(70,27)*{}; (70,33)*{} **\dir{--};
(80,27)*{}; (80,43)*{} **\dir{--};
(20,20)*{};
(20,30)*{\cdot} **\dir{--};
(20,40)*{\cdot} **\dir{--};
(20,50)*{\cdot} **\dir{--};
(20,60)*{\cdot} **\dir{--};
(20,70)*{} **\dir{--};
(90,20)*{};
(90,30)*{\cdot} **\dir{--};
(90,40)*{\cdot} **\dir{--};
(90,50)*{\cdot} **\dir{--};
(90,60)*{\cdot} **\dir{--};
(90,70)*{} **\dir{--};
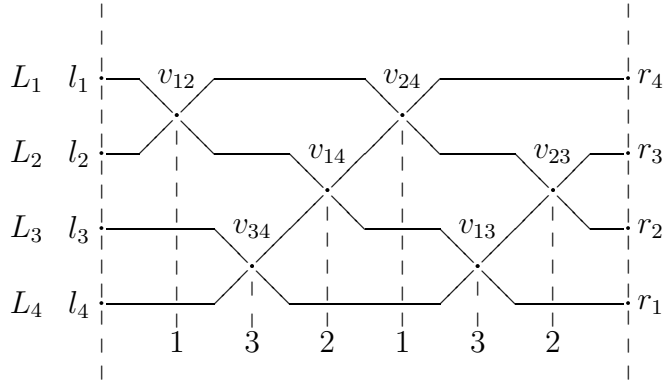
(50,15)*{\mbox{Figure 4.1 : wiring diagram of 
$\mathbf{w}_0=s_1 s_3 s_2 s_1 s_3 s_2$ of type $A_3$.}};
\end{xy}
\]
   
Let $G^{\circ}(\mathbf{w}_0)$ be the non-oriented graph
obtained from $\mathcal{WD}(\mathbf{w}_0)$, 
whose vertices are the crossing points of pseudolines, 
and whose edges are given by pseudoline
segments linking two crossing points. Likewise, 
let us denote by $G(\mathbf{w}_0)$ the non-oriented graph
obtained in a similar way, by considering 
the vertices of $G^{\circ}(\mathbf{w}_0)$, as well
as the vertices $l_1, \ldots , \, l_{n+1}, 
\, r_1, \ldots r_{n+1}$ on the border of 
$\mathcal{WD}(\mathbf{w}_0)$. 

Two vertices $v_{\beta_k}$ and $v_{\beta_{k'}}
$ with $k < k'$, 
$i_k=i$ and $i_{k'}=j$ are adjacent 
in $G^{\circ}(\mathbf{w}_0)$ if one 
of two cases occur. The adjacency
is \textbf{diagonal} if  
one has $\mid i-j \mid=1$ and
$i_l \neq i, \, i_l \neq j$ for any $l$
satisfying $k<l<k'$.
The adjacency is \textbf{horizontal}
 when $i=j$ and either $i_l < i-1$ for all $k < l < k'$, 
or $i_l > i+1$ for all $k < l < k'$. 
We see in Figure 4.1 above that  
$v_{12}$ and $v_{14}$ are diagonally adjacent, 
where as $v_{12}$ and $v_{24}$
are horizontally adjacent.

The wiring diagram $\mathcal{WD}(\mathbf{w}_0)$ defines a set 
of bounded chambers. Such a chamber may 
be indexed by set of indices of the pseudo-lines passing 
above it. This set may be seen as indices of a 
column Young tableau. One gets a one-to-one correspondence
between these chambers and a set of weights inside 
$\underset{k=1}{\overset{n}{\bigcup}} W \omega_k$.

A vertex $v_{\beta}$ of  $G^{\circ}(\mathbf{w}_0)$
may be assigned either the weight $\lambda^-(v_{\beta})$
of the chamber left to it,
or the weight $\lambda^+(v_{\beta})$  of the 
chamber right to it. 
It is well known that if $\beta=s_{i_1} s_{i_2} \ldots s_{i_{k-1}} (\alpha_{i_k})$,
then 
$$ \lambda^-(\beta)=s_{i_1} s_{i_2} \ldots s_{i_{k-1}} (\omega_{i_k}),
\; \; \; 
\lambda^+(\beta)=s_{i_1} s_{i_2} \ldots s_{i_{k-1}} s_{i_k} 
(\omega_{i_k}).$$

Figure 4.2 below gives the example of the chamber system for
$\mathbf{w}_0=s_1 s_3 s_2 s_1 s_3 s_2$ of type $A_3$.
 The vertex $v_{14}$ is between
the chamber labelled 12 to the left, and that
labelled 24 to the right. One has 
$\lambda^-(v_{14})=\nu_1+\nu_2, 
\; \lambda^+(v_{14})=\nu_2+\nu_4$.

\[
\begin{xy} 
(20,60)*{\cdot}; 
(25,60)*{} **\dir{-};
(30,55)*{\cdot} **\dir{-};
(35,50)*{} **\dir{-};
(45,50)*{} **\dir{-};
(50,45)*{\cdot} **\dir{-};
(55,40)*{} **\dir{-};
(65,40)*{} **\dir{-};
(70,35)*{\cdot} **\dir{-};
(75,30)*{} **\dir{-};
(85,30)*{} **\dir{-};
(90,30)*{\cdot} **\dir{-};
(20,50)*{\cdot};
(25,50)*{} **\dir{-};
(30,55)*{\cdot} **\dir{-};
(35,60)*{} **\dir{-};
(55,60)*{} **\dir{-};
(60,55)*{\cdot} **\dir{-};
(65,50)*{} **\dir{-};
(75,50)*{} **\dir{-};
(80,45)*{\cdot} **\dir{-};
(85,40)*{} **\dir{-};
(90,40)*{\cdot} **\dir{-};
(20,40)*{\cdot}; 
(35,40)*{} **\dir{-};
(40,35)*{\cdot} **\dir{-};
(45,30)*{} **\dir{-};
(65,30)*{} **\dir{-};
(70,35)*{\cdot} **\dir{-};
(75,40)*{} **\dir{-};
(80,45)*{\cdot} **\dir{-};
(85,50)*{} **\dir{-};
(90,50)*{\cdot} **\dir{-};
(20,30)*{\cdot};
(35,30)*{} **\dir{-};
(40,35)*{\cdot} **\dir{-};
(45,40)*{} **\dir{-};
(50,45)*{\cdot} **\dir{-};
(55,50)*{} **\dir{-};
(60,55)*{\cdot} **\dir{-};
(65,60)*{} **\dir{-};
(90,60)*{\cdot} **\dir{-};
(30,55)*{\cdot}; (40,35)*{\cdot}; 
(50,45)*{\cdot}; (60,55)*{\cdot}; 
(70,35)*{\cdot}; (80,45)*{\cdot};
(10,60)*{L_1}; (10,50)*{L_2}; (10,40)*{L_3}; (10,30)*{L_4};
(17,60)*{l_1}; (17,50)*{l_2}; (17,40)*{l_3}; (17,30)*{l_4};
(93,60)*{r_4}; (93,50)*{r_3}; (93,40)*{r_2}; (93,30)*{r_1};
(20,20)*{};
(20,30)*{\cdot} **\dir{--};
(20,40)*{\cdot} **\dir{--};
(20,50)*{\cdot} **\dir{--};
(20,60)*{\cdot} **\dir{--};
(20,70)*{} **\dir{--};
(90,20)*{};
(90,30)*{\cdot} **\dir{--};
(90,40)*{\cdot} **\dir{--};
(90,50)*{\cdot} **\dir{--};
(90,60)*{\cdot} **\dir{--};
(90,70)*{} **\dir{--};
(25,55)*{1};
(45,55)*{2};
(80,55)*{4};
(30,45)*{12};
(65,45)*{24};
(86,45)*{34};
(25,35)*{123};
(55,35)*{124};
(80,35)*{234};
\end{xy}
\]
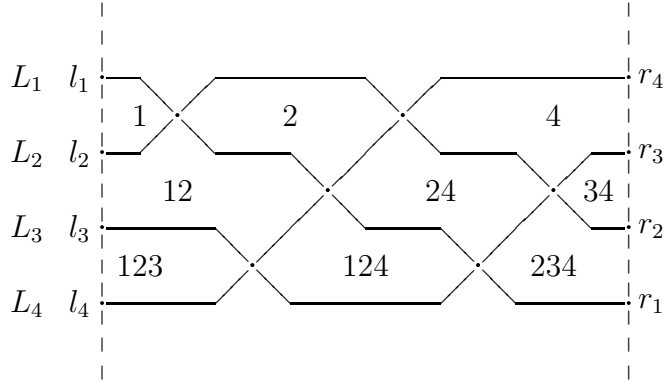
\hspace*{10mm} Figure 4.2 : chamber system of $\mathcal{WD}(\mathbf{w}_0)$,
for $\mathbf{w}_0=s_1 s_3 s_2 s_1 s_3 s_2$ of type $A_3$. \bigskip 

 Gleizer-Postnikov  (\cite{GleizerPostnikov} section 5) 
obtain a system of defining inequalities for the string cone 
$\mathcal{C}_{\mathbf{w}_0}$, 
by transforming for every $i \in I$, the graph $G(\mathbf{w}_0)$ 
into  an oriented graph $G(\wbold_0,i)$.
The pseudolines $L_1, \, L_2, \, \dots \, L_i$ are oriented backwards, and
the pseudolines $L_{i+1}, \ldots L_{n+1}$ forwards. 
The resulting graph is acyclic.

\[
\begin{xy} 
(20,60)*{\cdot}; 
(25,60)*{} **\dir{-} ?(0.5)*@{<};
(25,60)*{};
(30,55)*{\cdot} **\dir{-};
(35,50)*{} **\dir{-};
(45,50)*{} **\dir{-} ?(0.5)*@{<};
(45,50)*{};
(50,45)*{\cdot} **\dir{-};
(55,40)*{} **\dir{-};
(65,40)*{} **\dir{-} ?(0.5)*@{<};
(65,40)*{};
(70,35)*{\cdot} **\dir{-};
(75,30)*{} **\dir{-};
(90,30)*{\cdot} **\dir{-} ?(0.5)*@{<};
(20,50)*{\cdot};
(25,50)*{} **\dir{-} ?(0.5)*@{<};
(25,50)*{};
(30,55)*{\cdot} **\dir{-};
(35,60)*{} **\dir{-};
(55,60)*{} **\dir{-} ?(0.5)*@{<};
(55,60)*{};
(60,55)*{\cdot} **\dir{-};
(65,50)*{} **\dir{-};
(75,50)*{} **\dir{-} ?(0.5)*@{<};
(75,50)*{};
(80,45)*{\cdot} **\dir{-};
(85,40)*{} **\dir{-};
(90,40)*{\cdot} **\dir{-} ?(0.5)*@{<};
(20,40)*{\cdot}; 
(35,40)*{} **\dir{-} ?(0.5)*@{>};
(35,40)*{};
(40,35)*{\cdot} **\dir{-};
(45,30)*{} **\dir{-};
(65,30)*{} **\dir{-} ?(0.5)*@{>};
(65,30)*{};
(70,35)*{\cdot} **\dir{-};
(80,45)*{\cdot} **\dir{-} ?(0.5)*@{>};
(80,45)*{\cdot};
(85,50)*{} **\dir{-};
(90,50)*{\cdot} **\dir{-} ?(0.5)*@{>};
(20,30)*{\cdot};
(35,30)*{} **\dir{-} ?(0.5)*@{>};
(35,30)*{};
(40,35)*{\cdot} **\dir{-};
(50,45)*{\cdot} **\dir{-} ?(0.5)*@{>};
(50,45)*{\cdot};
(60,55)*{\cdot} **\dir{-} ?(0.5)*@{>};
(60,55)*{\cdot};
(65,60)*{} **\dir{-};
(90,60)*{\cdot} **\dir{-} ?(0.5)*@{>};
(30,55)*{\cdot}; (40,35)*{\cdot}; 
(50,45)*{\cdot}; (60,55)*{\cdot}; 
(70,35)*{\cdot}; (80,45)*{\cdot};
(10,60)*{L_1}; (10,50)*{L_2}; (10,40)*{L_3}; (10,30)*{L_4};
(17,60)*{l_1}; (17,50)*{l_2}; (17,40)*{l_3}; (17,30)*{l_4};
(93,60)*{r_4}; (93,50)*{r_3}; (93,40)*{r_2}; (93,30)*{r_1};
(20,20)*{};
(20,30)*{\cdot} **\dir{--};
(20,40)*{\cdot} **\dir{--};
(20,50)*{\cdot} **\dir{--};
(20,60)*{\cdot} **\dir{--};
(20,70)*{} **\dir{--};
(90,20)*{};
(90,30)*{\cdot} **\dir{--};
(90,40)*{\cdot} **\dir{--};
(90,50)*{\cdot} **\dir{--};
(90,60)*{\cdot} **\dir{--};
(90,70)*{} **\dir{--};
(30,60)*{v_{12}}; (40,40)*{v_{34}}; (50,50)*{v_{14}};
(60,60)*{v_{24}}; (70,40)*{v_{13}}; (80,50)*{v_{23}};
(50,15)*{\mbox{Figure 4.3 : 
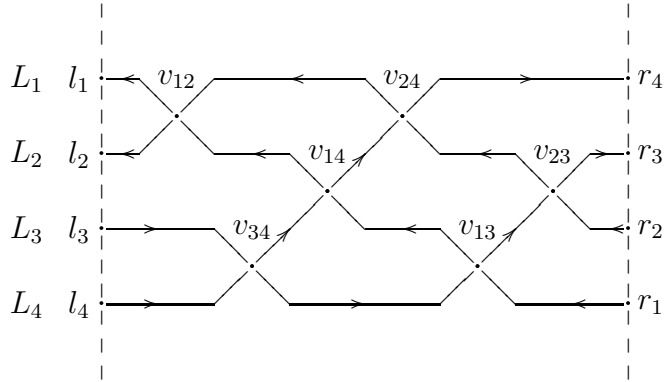
$G(\mathbf{w}_0,2), \; \mathbf{w}_0=s_1 s_3 s_2 s_1 s_3 s_2$, 
type $A_3$}};
\end{xy}
\]

Let $\pi$ be a path inside $G(\mathbf{w}_0,i)$. We shall
qualify the two configurations below as forbidden crossings :

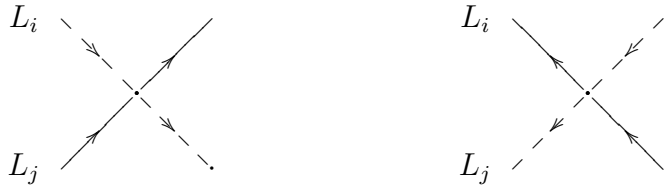
\[
\begin{xy}
(20,30)*{};
(30,40)*{\cdot} **\dir{-} ?(0.5)*@{>};
(30,40)*{\cdot};
(40,50)*{} **\dir{-} ?(0.5)*@{>};
(20,50)*{};
(30,40)*{\cdot} **\dir{--} ?(0.5)*@{>};
(30,40)*{\cdot};
(40,30)*{\cdot} **\dir{--} ?(0.5)*@{>};
(15,50)*{L_i};
(15,30)*{L_j};
(100,30)*{};
(90,40)*{\cdot} **\dir{-} ?(0.5)*@{>};
(90,40)*{\cdot};
(80,50)*{} **\dir{-} ?(0.5)*@{>};
(100,50)*{};
(90,40)*{\cdot} **\dir{--} ?(0.5)*@{>};
(90,40)*{\cdot};
(80,30)*{} **\dir{--} ?(0.5)*@{>};
(75,50)*{L_i};
(75,30)*{L_j};
(50,20)*{\mbox{(path $\pi$ in plain line)}};
\end{xy}
\]

A \textbf{Gleizer-Postnikov path (GP-path) of type} $i$
is a path of $G(\mathbf{w}_0,i)$, starting from
$l_{i+1}$ on the left border, ending at $l_i$, and
that does not contain
forbidden crossings.

An example of such a path in Figure 4.3 above is 
$l_3 \longrightarrow v_{34} \longrightarrow v_{13}
     \longrightarrow v_{14} \longrightarrow v_{24}
     \longrightarrow v_{12} \longrightarrow l_2$.  
There are $5$ such paths inside $G(\mathbf{w}_0,2)$. 

\textbf{Remark : } Gleizer-Postnikov paths were
called rigorous paths in \cite{GleizerPostnikov}.
We have translated the vertical setting of 
\cite{GleizerPostnikov} 
to an horizontal one, which is
more natural when comparing wiring diagrams
to Auslander-Reiten quivers.
 
Let $\pi$ be a $GP$-path. 
If $\pi$ enters a vertex $v_{\beta_j}$ following 
the line $L_h$, and leaves it following $L_l$, 
assign to it the value
$$ k_j:=\left\{ 
\begin{array}{rl} 
1 & \mbox{if} \;  h>l \\
-1 & \mbox{if} \; h<l \\
0 & \mbox{if} \;  h=l  \\
\end{array}
\right.$$ 
If $v_{\beta_j}$ is not a vertex of $\pi$, put $k_j:=0$. 
The coordinates $k_j, \, j=1, \ldots, \, N$ define
a vector $\kbold_{\pi}$ of $\Zbold^N$.

Take as an example the GP-path above,
$\pi=l_3 \longrightarrow v_{34} \longrightarrow v_{13}
     \longrightarrow v_{14} \longrightarrow v_{24}
     \longrightarrow v_{12} \longrightarrow l_2$.
The path $\pi$ starts on strand
$L_3$, changes to strand $L_1$ at $v_{13}$ (hence
a positive contribution), then passes from
strand $L_1$ to strand $L_4$ at $v_{14}$ (hence
a negative contribution), and finally
its last change of strands occurs at $v_{24}$,
where $\pi$ passes from $L_4$ to $L_2$ 
(hence a positive contribution). The
vertices $v_{13}, \, v_{14}, \, v_{24}$
occur respectively in positions 5, 3  and 4 
while we follow $\mathcal{WD}(s_1 s_3 s_2 s_1 s_3 s_2)$ 
from left to right. We get 
$\mathbf{k}_{\pi}=(0,0,-1,1,1,0)$.

As $G(\wbold_0,i)$ is acyclic, there are only
finitely many different Gleizer-Postnikov paths
of type $i$.

\textbf{Theorem 4.1} (\cite{GleizerPostnikov} corollary 5.8)

\emph{
Let $K^{GP}_{\mathbf{w}_0}$ be the set
of all vectors $\kbold_{\pi}$, where
$\pi$ varies over
all possible GP-paths of all possible types 
$i \in I$. Then $K^{GP}_{\mathbf{w}_0}$
defines the string cone $\mathcal{C}_{\mathbf{w}_0}$. 
One has 
$$ \mathcal{S}_{\mathbf{w}_0}=\{ \tbold \in \Nbold^N \mid
 \forall \kbold_{\pi} \in K^{GP}_{\mathbf{w}_0}, \; \; 
\kbold_{\pi} \cdot \tbold \geq 0 \}.$$}

Fix $i \in I$, and consider $G(\wbold_0,i)$.
Let us denote by $\delta_i^>$ the 
segment of $L_{i+1}$ starting from 
the left border on $l_{i+1}$, 
up to its intersection   
$v_{\alpha_i}$ with $L_i$. In a similar way,
we denote by $\delta_i^<$ the segment of $L_i$
starting from $v_{\alpha_i}$ and going
back to $l_i$ on the left border following $L_i$. 
The concatenation $\delta_i:=\delta_i^>*\delta_i^<$ is
then a path starting from $l_{i+1}$ and ending
in $l_i$. We shall call it the \textbf{limiting path of type} $i$.
 
Observe $\delta_i$ serves as the boundary
of a set of chambers $\mathcal{Z}_i(\mathbf{w}_0)$ .
These chambers are those lying below $L_i$, 
above $L_{i+1}$ and left of $v_{\alpha_i}$. 
Let us denote by $Z_i(\mathbf{w}_0)$ 
the set of vertices of $G^{\circ}(\mathbf{w}_0)$
which are the rightmost vertices of chambers
of $\mathcal{Z}_i(\mathbf {w}_0)$. 

\textbf{Lemma 4.2}

\emph{
Consider $v$ a vertex of $\delta_i$ other than 
$l_i, \, l_{i+1}, \, v_{\alpha_i}$. Suppose this
vertex is a crossing of $\delta_i$ with a line 
$L_k \, (k \neq i, i+1)$. We have then the
two following cases :
\begin{itemize}
\item[i)] If $v \in \delta_i^>$, then $L_k$ 
crosses $\delta_i$ going inside $\mathcal{Z}_i(\mathbf{w}_0)$.
\item[ii)] If $v \in \delta_i^<$, then $L_k$
crosses $\delta_i$ going outside $\mathcal{Z}_i(\mathbf{w}_0)$.
\end{itemize}}

\vfill \eject

\textbf{Proof of lemma 4.2 :}

i) We have to eliminate the possibility of $L_k$
going out of $\mathcal{Z}_i(\mathbf{w}_0)$ while
crossing $\delta_i^>$. By definition, $\mathcal{Z}_i(\mathbf{w}_0)$ lies
above $\delta_i^>$, hence
$L_k$ crosses $\delta_i^>$ going downwards.
Depending on the orientation of $L_k$, we
obtain the two following cases :
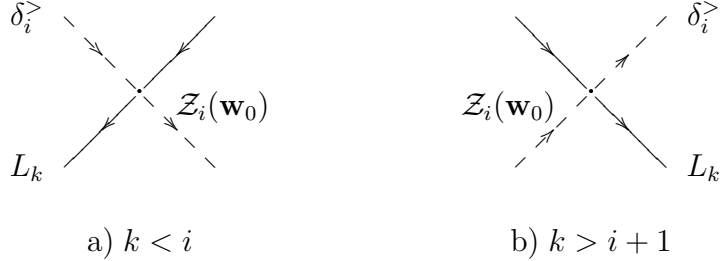
\[
\begin{xy}
(40,50)*{};
(30,40)*{\cdot} **\dir{-} ?(0.5)*@{>};
(30,40)*{\cdot};
(20,30)*{} **\dir{-} ?(0.5)*@{>};
(20,50)*{};
(30,40)*{\cdot} **\dir{--} ?(0.5)*@{>};
(30,40)*{\cdot};
(40,30)*{} **\dir{--} ?(0.5)*@{>};
(15,30)*{L_k};
(15,50)*{\delta_i^>};
(41,38)*{\mathcal{Z}_i(\mathbf{w}_0)};
(30,20)*{\mbox{a)} \; k<i}; 
(80,50)*{};
(90,40)*{\cdot} **\dir{-} ?(0.5)*@{>};
(90,40)*{\cdot};
(100,30)*{} **\dir{-} ?(0.5)*@{>};
(80,30)*{};
(90,40)*{\cdot} **\dir{--} ?(0.5)*@{>};
(90,40)*{\cdot};
(100,50)*{} **\dir{--} ?(0.5)*@{>};
(105,50)*{\delta_i^>};
(105,30)*{L_k};
(79,38)*{\mathcal{Z}_i(\mathbf{w}_0)};
(90,20)*{\mbox{b)} \; k>i+1};
\end{xy}
\]

In case a), $L_k$ would have to recross $L_{i+1}$
in order to return to  $l_k$ which lies above 
$l_{i+1}$ on the left border. In case b), $l_k$
lies below $l_{i+1}$ on the left border, so 
$L_k$ would have to cross $L_{i+1}$ a first time
in order to reach the crossing point in $b)$ from above. 
Both cases are impossible since they would imply
at least two crossings of $L_k$ and $L_{i+1}$
in $\mathcal{WD}(\mathbf{w}_0)$.

Statement ii) is proved by symmetric arguments. $\Box$

\textbf{Corollary 4.3}

\emph{
$\delta_i$ is a $GP$-path.}

\textbf{Proof :}

Clearly, $v_{\alpha_i}$ is not a forbidden crossing.
The other vertices of $\delta_i$ belong to the cases
detailed in the lemma above, none of
them being forbidden. $\Box$
\medskip

\textbf{Proposition 4.4}

\emph{
Let $\pi$ be a GP-path of type $i$. Then $\pi$ stays
inside $\mathcal{Z}_i(\mathbf{w}_0)$.}

\textbf{Proof : }

The path $\pi$ starts at $l_{i+1}$ and ends at $l_i$ which
are inside $\mathcal{Z}_i(\mathbf{w}_0)$. Suppose it exits 
$\mathcal{Z}_i(\mathbf{w}_0)$ at some
vertex $v_1 \in \delta_i$. By the lemma above, one
must have $v_1 \in \delta_i^<$. The same lemma shows
$\pi$ must return inside $\mathcal{Z}_i(\mathbf{w}_0)$ through a crossing point
$v_2 \in \delta_i^>$. Now the segment of $\delta_i$
between $v_1$ and $v_2$ goes from $v_2$ to $v_1$, so
we can create a cycle. This is
in contradiction with the fact that $G(\wbold_0,i)$
is acyclic. $\Box$

\section*{5. Hammocks}

\hspace*{5mm} We continue to restrict ourselves to $A_n$ case, with
$Q$ and $\mathbf{w}_0$ fixed. 
Let us denote  by $\lbrack \beta : \alpha_i \rbrack$ 
the coefficient of $\alpha_i$ in the expression of $\beta$.
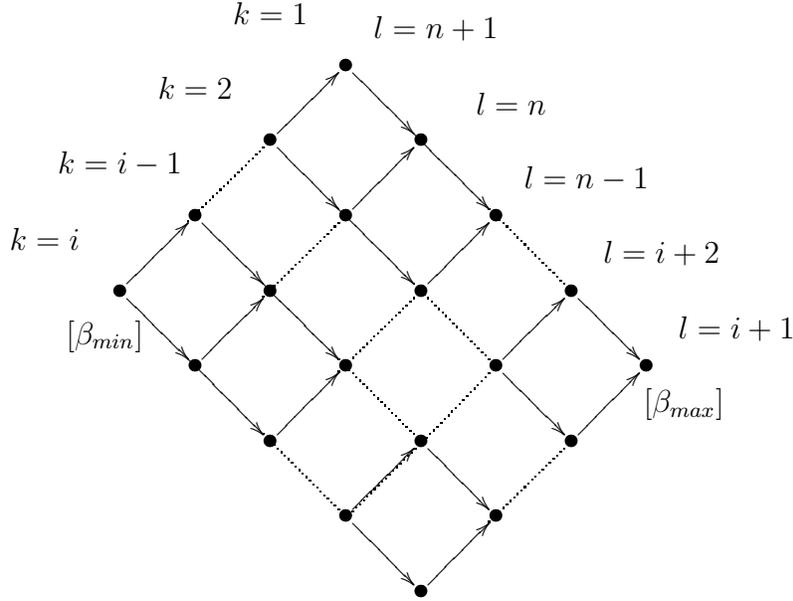
The \textbf{hammock} of type $i$, $i \in I$   
(\cite{Brenner}) is 
$H_i(Q):=\{ \lbrack \beta \rbrack \mid 
\lbrack \beta : \alpha_i \rbrack >0 \}$. 
The set $P_i(Q)$ is a subset of $H_i(Q)$, and
we shall see the combinatorics
of Lusztig's moves of type $i$ is obtained from that
of $H_i(Q)$. The structure of $H_i(Q)$ itself is very
simple, and is deduced from the Coxeter element $c$
attached to $Q$.  

Recall \cite{BernsteinGelfandPonomarev} that
one may renumber the vertices of $Q$ by a 
permutation $i_1, \, i_2, \ldots, \, i_n$  
of $I$, such that for every arrow $j \longrightarrow k$
of $Q$, one has $i_j < i_k$. 
The Coxeter element is then given by $c=s_{i_1} s_{i_2} \ldots s_{i_n}$.
 The action of $c$ on
$\Phi^+$ is the mirror image of the action of
the translation $\tau$ upon $\Gamma_Q$ : 
if $N=\tau M$ then $\dbold_N=c \dbold_M$.

In type $A_n$, as $W \cong \gothicS_{n+1}$, $c$ is an 
$n+1$-cycle. Its expression may be constructed 
out of $Q$ by the following algorithm 
(\cite{Zelikson} Lemma 4.2) :
\begin{itemize}
\item[-] Start with just the element $n+1$.
\item[-] Proceed in decreasing order 
$i=n, n-1, \ldots, 2$ :
\begin{itemize}
\item[-] If inside $Q$, one has 
$\stackrel{i-1}{\cdot} \longleftarrow \stackrel{i}{\cdot}$,
add $i$ to the right of the indices already written.
\item[-] If inside $Q$, one has 
$\stackrel{i-1}{\cdot} \longrightarrow \stackrel{i}{\cdot}$,
add $i$ to the left of the indices already written.
\end{itemize}
\item[-] Finish by adding $1$ to the left of the $n$ indices already
written.
\end{itemize}  

The sequence of indices $j_1 \, j_2 \, \ldots j_{n+1}$ thus
obtained is a $n+1$-cycle expression of $c$. The consequence of this specific
algorithm is that the $n+1$-cycle expressions of $c$ 
verify a special "segment" property.
\medskip

\textbf{Lemma 5.1} (\cite{Zelikson}, Proposition 4.3) 

\emph{
For every  $i \in \{ 1, \, \ldots, \,  n \}$, there is
a cycle expression $c=(j_1, \ldots , j_{i}, j_{i+1} 
\ldots j_{n+1})$ where : 
\begin{itemize}
\item[-]  $j_1, j_2 , \ldots j_{i}$ is a permutation
of $1, \ldots i$.
\item[-] $j_{i+1}, j_{i+2}, \ldots j_{n+1}$ 
is a permutation of $i+1, i+2, \ldots n+1$.
\end{itemize}}

We shall refer to this  cycle expression 
as the $i$-\textbf{segmented expression} of $c$, and 
denote it by 
$(j_1 , \ldots j_i \mid j_{i+1}, \ldots , j_{n+1})$.
\medskip

\textbf{Example}

Consider the quiver 
$Q \; :  \; 
     \stackrel{1}{\cdot} \longleftarrow
     \stackrel{2}{\cdot} \longrightarrow
     \stackrel{3}{\cdot} \longleftarrow
     \stackrel{4}{\cdot}$ of type $A_4$. 
One has $c=s_2 s_1 s_4 s_3$. The algorithm
above gives 
$$  5 \stackrel{right}{\longrightarrow} 
   54 \stackrel{left}{\longrightarrow}
  354 \stackrel{right}{\longrightarrow} 
 3542 \stackrel{left}{\longrightarrow}
13542.
  $$
One may verify the 5-cycle expression $(13542)$ 
obtained agrees with $c$.  
The $i$-segmented expressions of $c$ are then respectively
$$ \begin{array}{cc}
i=1 \; : \;  (1 \mid 3542) &  \hspace*{5mm} i=3 \; : \;  (213 \mid 54) \\
i=2 \; : \;  (21 \mid 354) &  \hspace*{5mm} i=4 \; : \;  (4213 \mid 5)
\end{array} $$

The positions of $i, i+1$
in the $i$-segmented expression of $c$ depend 
on the \linebreak[4] neighbourhood of $i$ in $Q$ : 
$$ \begin{array}{rlrl}
\stackrel{i-1}{\cdot} \longrightarrow \stackrel{i}{\cdot} \; 
\mbox{in Q} \; : & j_i=i & 
\hspace*{5mm} \stackrel{i-1}{\cdot} \longleftarrow \stackrel{i}{\cdot} \;
\mbox{in Q} \; : & j_1=i \\
\stackrel{i}{\cdot} \longrightarrow \stackrel{i+1}{\cdot} \;
\mbox{in Q} \; : & j_{i+1}=i+1 &
\hspace*{5mm} \stackrel{i}{\cdot} \longleftarrow \stackrel{i+1}{\cdot} \;
\mbox{in Q} \; : & j_{n+1}=i+1 
\end{array} $$
In  the example above, the neighbourhood of vertex 2 in $Q$ is
$\stackrel{1}{\cdot} \longleftarrow \stackrel{2}{\cdot} \longrightarrow 
\stackrel{3}{\cdot}$, so that $2$ appears in first position of the segment $21$, and
$3$ appears in first position of the segment $354$.

\vfill \eject

Given an integer $m \geq 1$, we shall denote by 
$\lbrack m \rbrack$ the set $\{ 1, \, 2, \, \ldots m \}$. 
\medskip

\textbf{Proposition 5.2}
\emph{
Let $Q$ be a quiver of type $A_n$.
\begin{itemize} 
\item[a)] The set $H_i(Q)$ for i $\in I$, seen as a subgraph of $\Gamma_Q$, 
has a structure isomorphic to \linebreak[4] 
$\lbrack i \rbrack \times \lbrack n+1-i \rbrack$ 
as given in Figure 5.1, 
with $\lbrack \beta_{min} \rbrack$ the isoclass
of the projective cover of the simple module $S_i$
and $\lbrack \beta_{max} \rbrack$ the isoclass
of the injective envelope of $S_i$.
\item[b)] If 
$(j_1 j_2 \ldots j_i \mid j_{i+1} \ldots j_{n+1})$ is
the $i$-segmented writing of $c$, then
the vertex at position $(k,l)$ in the Figure 5.1
is $\lbrack \alpha_{j_k}+\alpha_{j_{k}+1}+ \ldots
\alpha_{j_l-1} \rbrack$.
\end{itemize}}

\[
\begin{xy}
(55,20)*{\bullet};
(65,30)*{\bullet} **\dir{-} ?>*@{>};
(75,40)*{\bullet} **\dir{.};
(85,50)*{\bullet} **\dir{-} ?>*@{>};
(90,45)*{\lbrack \beta_{max} \rbrack};
(97,55)*{l=i+1};
(45,30)*{\bullet};
(55,40)*{\bullet} **\dir{-} ?>*@{>};;
(65,50)*{\bullet} **\dir{.};
(75,60)*{\bullet} **\dir{-} ?>*@{>};;
(87,65)*{l=i+2};
(35,40)*{\bullet};
(45,50)*{\bullet} **\dir{-} ?>*@{>};
(55,60)*{\bullet} **\dir{.};
(65,70)*{\bullet} **\dir{-} ?>*@{>};
(77,75)*{l=n-1};
(25,50)*{\bullet};
(35,60)*{\bullet} **\dir{-} ?>*@{>};
(45,70)*{\bullet} **\dir{.};
(55,80)*{\bullet} **\dir{-} ?>*@{>};
(67,85)*{l=n};
(15,60)*{\bullet};
(25,70)*{\bullet} **\dir{-} ?>*@{>};
(35,80)*{\bullet} **\dir{.};
(45,90)*{\bullet} **\dir{-} ?>*@{>};
(57,95)*{l=n+1};
(5,67)*{k=i};
(13,54)*{\lbrack \beta_{min} \rbrack};
(15,60)*{\bullet};
(25,50)*{\bullet} **\dir{-} ?>*@{>};
(35,40)*{\bullet} **\dir{-} ?>*@{>};
(45,30)*{\bullet} **\dir{.};
(55,20)*{\bullet} **\dir{-} ?>*@{>};
(15,77)*{k=i-1};
(25,70)*{\bullet};
(35,60)*{\bullet} **\dir{-} ?>*@{>};
(45,50)*{\bullet} **\dir{-} ?>*@{>};
(55,40)*{\bullet} **\dir{.};
(65,30)*{\bullet} **\dir{-} ?>*@{>};
(25,87)*{k=2};
(35,80)*{\bullet};
(45,70)*{\bullet} **\dir{-} ?>*@{>};
(55,60)*{\bullet} **\dir{-} ?>*@{>};
(65,50)*{\bullet} **\dir{.};
(75,40)*{\bullet} **\dir{-} ?>*@{>};
(35,97)*{k=1};
(45,90)*{\bullet};
(55,80)*{\bullet} **\dir{-} ?>*@{>};
(65,70)*{\bullet} **\dir{-} ?>*@{>};
(75,60)*{\bullet} **\dir{.};
(85,50)*{\bullet} **\dir{-} ?>*@{>};
(50,10)*{\mbox{Figure 5.1 : Hammock $H_i(Q)$ of 
type $A_n$}};
\end{xy}
\]

The proposition consists in a computation well
known to specialists. The reader may consult 
\cite{GabrielARquivers} section 6.5 for technical details,
especially pages 52-54 (our Figure 5.1 corresponds
to the first scheme in \cite{GabrielARquivers} Figure 15).
We provide here some guidelines for non-specialists.

The set $H_i(Q)$ is particularly simple to compute in the case
of the "$i$-regular" quiver
$ \stackrel{1}{\cdot} \longrightarrow 
   \stackrel{2}{\cdot} \longrightarrow 
   \ldots 
   \longrightarrow 
   \stackrel{i}{\cdot} 
   \longleftarrow
   \ldots
   \longleftarrow
   \stackrel{n-1}{\cdot} \longleftarrow
   \stackrel{n}{\cdot}$
 admitting $i$
as its unique sink. 
The modules at positions $(1,n+1), \, (2,n+1), \, \ldots, \, (i,n+1), \, (i,n), \, 
\ldots , \, (i,i+1)$ in Figure 5.1 are the respective projective covers of the simple modules
$S_1, \, \ldots, \,  S_n$. Their dimension vectors are directly obtained out of $Q$.

The mesh structure of $\Gamma_Q$ \cite{GabrielARquivers},
verifies additivity for dimension vectors. For each "square" configuration
of vertices $\lbrack X \rbrack, \, \lbrack{Z_1} \rbrack , \, \lbrack{Z_2} \rbrack, \, 
\lbrack Y \rbrack$,
respectively at positions \linebreak[4] $(k,l), \, (k-1,l), \, (k,l-1), \, (k-1,l-1)$ in Figure 5.1
one has $\dbold_X+\dbold_Y=\dbold_{Z_1}+\dbold_{Z_2}$. 
Thus the knowledge of the dimension vectors for the "slice" of projective modules 
allows to compute the rest of $H_i(Q)$. 

An arbitrary quiver $Q$ may be obtained from
the $i$-regular one by a sequence of orientations
changes transforming a source vertex $j$, with $j \neq i$, 
into a sink. The effect of this transformation 
on vertices of $\Gamma_Q$ other than $\lbrack \alpha_j \rbrack$
consists in applying the corresponding BGP-reflection functor 
$\Sigma_j$ \cite{BernsteinGelfandPonomarev}. 
This is the case of the vertices of $H_i(Q)$, which never
contains $\lbrack \alpha_j \rbrack$. 
In terms of dimension vectors, $\dbold_{\Sigma_j X}=s_j \dbold_X$.
One checks this change is coherent with the change in the Coxeter element $c$
 due to the change of orientation.

One may verify by direct computation that $P_i(Q)$ is
the order ideal defined by $\lbrack \alpha_i \rbrack$
inside $H_i(Q)$. The possible cases,  
according to the neighbourhood of $i$ inside $Q$, are :
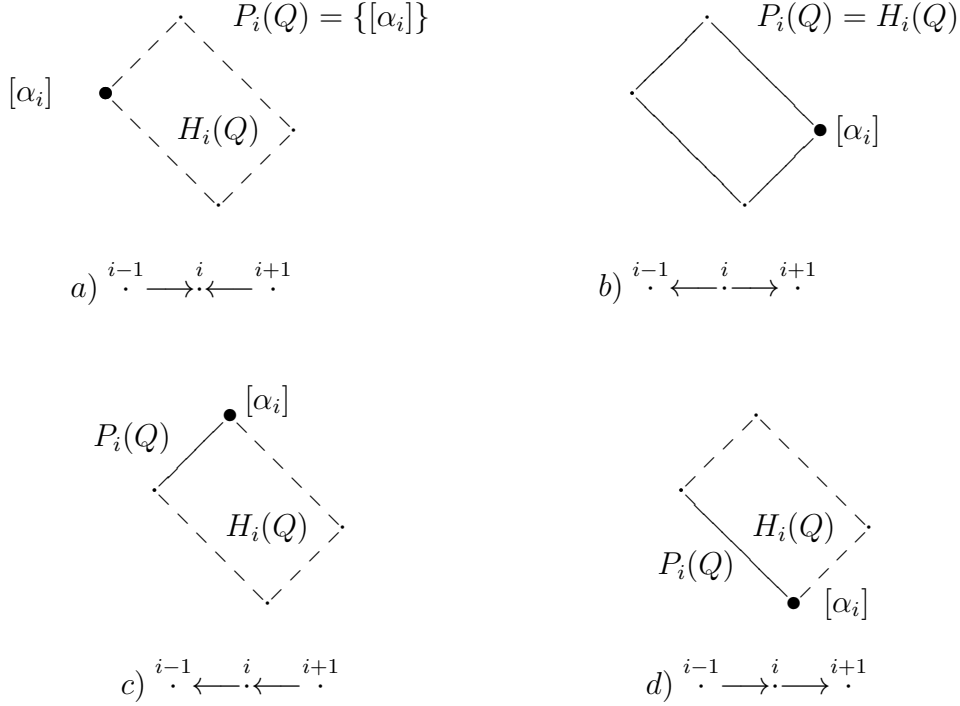
\[
\begin{xy}
(5,80)*{\lbrack \alpha_i \rbrack};
(15,80)*{\bullet};
(25,90)*{\cdot} **\dir{--};
(40,75)*{\cdot} **\dir{--};
(30,65)*{\cdot} **\dir{--};
(15,80)*{\bullet} **\dir{--};
(30,75)*{H_i(Q)};
(45,90)*{P_i(Q)=\{ \lbrack \alpha_i \rbrack \}};
(25,55)*{a)\stackrel{i-1}{\cdot} \longrightarrow
           \stackrel{i}{\cdot} \longleftarrow 
           \stackrel{i+1}{\cdot}};
(85,80)*{\cdot};
(95,90)*{\cdot} **\dir{-};
(110,75)*{\bullet} **\dir{-};
(100,65)*{\cdot} **\dir{-};
(85,80)*{\cdot} **\dir{-};
(115,90)*{P_i(Q)=H_i(Q)};
(115,75)*{\lbrack \alpha_i \rbrack};
(95,55)*{b)\stackrel{i-1}{\cdot} \longleftarrow
           \stackrel{i}{\cdot} \longrightarrow 
           \stackrel{i+1}{\cdot}};
\end{xy}
\]

\[
\begin{xy}
(15,30)*{\cdot};
(25,40)*{\bullet} **\dir{-};
(40,25)*{\cdot} **\dir{--};
(30,15)*{\cdot} **\dir{--};
(15,30)*{\cdot} **\dir{--};
(12,37)*{P_i(Q)};
(30,42)*{\lbrack \alpha_i \rbrack};
(30,25)*{H_i(Q)};
(25,5)*{c)\stackrel{i-1}{\cdot} \longleftarrow
           \stackrel{i}{\cdot} \longleftarrow 
           \stackrel{i+1}{\cdot}};
(85,30)*{\cdot};
(95,40)*{\cdot} **\dir{--};
(110,25)*{\cdot} **\dir{--};
(100,15)*{\bullet} **\dir{--};
(85,30)*{\cdot} **\dir{-};
(87,20)*{P_i(Q)};
(107,15)*{\lbrack \alpha_i \rbrack};
(100,25)*{H_i(Q)};
(95,5)*{d)\stackrel{i-1}{\cdot} \longrightarrow
           \stackrel{i}{\cdot} \longrightarrow 
           \stackrel{i+1}{\cdot}};
\end{xy}
\]

The correspondence $\Psi \, : \, \Gamma_Q \longrightarrow G^{\circ}(\mathbf{w}_0)$,
$\lbrack \beta \rbrack \mapsto v_{\beta}$
will allow us to transfer the combinatorial results
above to the setting of 
$\mathcal{WD}(\mathbf{w}_0)$. In particular,
the translation operation $\tau$ on $\Gamma_Q$
corresponds under $\Psi$, to an operation
$\tau_{wd}$ on vertices of $G^{\circ}(\mathbf{w}_0)$.
It is well known (\cite{Bedard}, Lemma 2.11), 
that for every $i \in I$, $\Psi$ establishes an
order preserving one-to-one correspondence 
between the $i$-th translation level of $\Gamma_Q$ and
crossings on the $i^{th}$ level of $\mathcal{WD}(\mathbf{w}_0)$.
Given $k$ between $1$ and $N$, let $k^-$ denote
the maximal index $j$ such that $j<k$ and $s_{i_j}=s_{i_k}$ in $\mathbf{w}_0$
(if such an index exists). The isoclass $\lbrack \beta_k \rbrack$ is not projective 
exactly when $k^{-}$ is defined. One has then 
$\tau_{wd}(v_{\beta_k})=v_{\beta_{k^{-}}}=v_{c \beta_k}$. 

Let us define $H_i(\mathbf{w}_0):=\Psi(H_i(Q))$.

\textbf{Proposition 5.3}

\emph{
Fix $i \in I$ and
consider $H_i(Q)$  and $H_i(\mathbf{w}_0)$ respectively 
as subgraphs
of $\Gamma_Q$ and $G^{\circ}(\mathbf{w}_0)$.
Then $\Psi$ restricted to $H_i(Q)$ establishes
an isomorphism of non-oriented graphs between
$H_i(Q)$  and $H_i(\mathbf{w}_0)$.}

\textbf{Proof : }

Let $\lbrack \beta_k \rbrack \longrightarrow 
\lbrack \beta_{k'} \rbrack, \; k<k'$ be an arrow of 
$H_i(Q)$. By \cite{Zelikson} Propostion 1.2, 
the vertices $i_k$ and $i_{k'}$ are linked in $D$,
and there is no 
occurrence either of a reflection $s_{i_j}=s_{i_k}$
or of a reflection $s_{i_j}=s_{i_{k'}}$
in the positions $j$ between $k+1$ and $k'-1$. 
This means that 
$v_{\beta_k}$ and $v_{\beta_{k'}}$ are
diagonally adjacent inside $H_i(\mathbf{w}_0)$.

Suppose now $v_{\beta_k}$ and $v_{\beta_{k'}}, \; k<k'$, are adjacent
vertices of $H_i(\mathbf{w}_0)$. A reduced expression $\mathbf{w}_0$
adapted to a quiver is alternating (\cite{Zelikson} Lemma 1.4), 
in the sense that if $j$ is linked to $l$ in $D$, then
between two occurrences 
of a reflection $s_j$, there is one and only one occurrence
of a reflection $s_l$. 
This excludes the possibility of
horizontal adjacencies inside
$H_i(\mathbf{w}_0)$.
We have  therefore 
$\beta_k=s_{i_1} \ldots s_{i_{k-1}} (\alpha_{i_k}), \;
\beta_{k'}=s_{i_1} \ldots s_{i_{k'-1}} (\alpha_{i_{k'}})$, 
with $(\alpha_{i_k}, \alpha_{i_{k'}})=-1$. 
One has
$(\beta_k, \beta_{k'})=
(\alpha_{i_k}, s_{i_k} s_{i_{k+1}} \ldots s_{i_{k'-1}} 
(\alpha_{i_{k'}}))$ by invariance of the Cartan scalar
product under action of $W$.
Now $i_{k+1}, \ldots i_{k'-1}$ are all different from
$i_k$ and $i_{k'}$, which implies that 
$(\alpha_{i_k}, s_{i_k} s_{i_{k+1}} \ldots s_{i_{k'-1}} 
(\alpha_{i_{k'}}))=1$. 

By \cite{Bedard} Lemma 2.11, $(\beta_k, \beta_{k'})=1$ means there is a 
path inside $\Gamma_Q$ from $\lbrack \beta_k \rbrack$
to $\lbrack \beta_{k'} \rbrack$. Now this path may 
have only one vertex on each of the adjacent levels
$i_k$ and $i_k'$, so it can cross them only once. 
In view of the mesh structure of $\Gamma_Q$, 
the only possibility is that this path is reduced 
to a single arrow $\lbrack \beta_k \rbrack \longrightarrow
\lbrack \beta_{k'} \rbrack$. $\Box$
\medskip

\textbf{Proposition 5.4}

\emph{
One has $Z_i(\mathbf{w}_0)=\Psi(P_i(Q))$. Furthermore,
these sets, seen as non-oriented \linebreak[4] subgraphs 
respectively of $\Gamma_Q$ and $G^{\circ}(\mathbf{w}_0)$,
are isomorphic.}

Only the first part of the proposition needs to 
be proved, the second then follows, using
Proposition 5.3. We shall need two results concerning
the matrix $R=(r_{i,j})$ of the Ringel form of $Q$.  
Let us define for $i\in I$, 
$\rho_i:=\underset{k=1}{\overset{n}{\sum}} r_{k,i} \omega_k$.
The weight $\rho_i$ is such that its coordinates
in the basis of fundamental weights are
given by the $i$-th column of $R$.
In a similar manner, let us put, corresponding 
to the $i$-th line of $R$, 
$\rho_i^t:=\underset{k=1}{\overset{n}{\sum}} r_{i,k} \omega_k$.

\textbf{Lemma 5.5} 
 
\emph{
The action of the Coxeter element on the vectors $\rho_i, \, i \in I$  
is given by $c \rho_i=-\rho_i^{t} $.}

\textbf{Proof :} 

Case by case analysis, following the $4$ cases
for the neighbourhood of $i$ in $Q$. Our convention
is $\omega_0=\omega_{n+1}=0$.

\underline{Case 1} : i is a source of $Q$

We have in this case
$\rho_i=\omega_i, \; 
\rho_i^{t} =-\omega_{i-1}+\omega_i-\omega_{i+1}$
and $c=s_i c_1$, $c_1$ being a 
product of simple reflections
$s_j$ with $j \neq i$. As $s_j(\omega_i)=\omega_i$ for
$j \neq i$, we see 
$c \rho_i=s_i (\omega_i)=
\omega_{i-1}-\omega_i+\omega_{i+1}$.

\underline{Case 2} : The neighbourhood of $i$ in $Q$ is of
the form  
$\stackrel{i-1}{\cdot} \longrightarrow
 \stackrel{i}{\cdot} \longrightarrow
 \stackrel{i+1}{\cdot}$

We have in this case 
$\rho_i=-\omega_{i-1}+\omega_i, \,
\rho_i^t=\omega_i-\omega_{i+1}$ and
$c=c_1 s_i c_2$, where $c_1$ is a product of 
simple reflections $s_j$ with $j \geq i+1$, and
$c_2$ a product of simple reflections $s_j$ with
$j \leq i-1$. We get 
$$ \begin{array}{rl}
 c_1 \rho_i &=\rho_i, \\
 s_i c_1 \rho_i &= s_i \rho_i=-\omega_i+\omega_{i+1}\\
 c_2 s_i c_1 \rho_i &=-\omega_i+\omega_{i+1}.
\end{array}$$

The other two cases follow by similar computations. $\Box$

\vfill  \eject

Recall we associated to a vertex 
$v_{\beta} \in G^{\circ }(\mathbf{w}_0)$, the weight
$\lambda^+(v_{\beta})$.
Let $\varphi_R$ be the linear mapping of the Euclidean
space $\Rbold \Phi$ defined by 
$\varphi_R (\alpha_i)=-\rho_i, \; i \in I$.

\textbf{Theorem 5.6} \emph{(\cite{Zelikson} Theorem 2.4)}

$$\forall \beta \in \Phi^+, \; \;
\lambda^+(v_{\beta})=\varphi_R(\beta).$$

The key of the theorem is \cite{Zelikson} Lemma 2.2 : if 
$k$ is a sink of an $ADE$ quiver $Q$, $Q':=s_k Q$, and 
$R, R'$ are the respective Ringel matrices, then
$s_k \varphi_{R'}=\varphi_R s_k$. This lemma implies
$\varphi_R$ commutes with $c^{-1}$, and hence
with $c=(c^{-1})^n$. If $\beta \in \Phi^+$ is such
that $c \beta \in \Phi^+$ (that is, if
$\tau$ is defined on $\lbrack \beta \rbrack$), 
then 
$c \lambda^+(v_{\beta})=c \varphi_R(\beta)=\varphi_R(c \beta)$. 
We have then that 
$c \lambda^+(v_{\beta})=\lambda^+(v_{c \beta})=\lambda^-(v_{\beta})$.

\textbf{Proof of Proposition 5.4}

For any $\beta \in \Phi^+$, one has that 
$(\beta, \alpha_i)_R =(\beta,\rho_i)$. This implies, 
by Theorem 3.1, 
that $\lbrack \beta \rbrack \in P_i(Q)$ if and only if 
$(\beta, \rho_i) >0$. 

Let us apply a similar analysis for a vertex $v_{\beta}$ 
of $Z_i(\mathbf{w}_0)$. 
A chamber is in $\mathcal{Z}_i(\mathbf{w}_0)$
if the strand $L_i$ passes above it, and the strand
$L_{i+1}$ passes below it. This means the index $i$
appears in the labelling of the chamber, and that
the index $i+1$ does not appear. 
Now the index $i$ corresponds to the weight 
$-\omega_{i-1}+\omega_{i}$ and
the index $i+1$ to the weight 
$-\omega_{i}+\omega_{i+1}$
The condition above in terms
of strands, translates in terms of weights, 
to having $\omega_i$ with coefficient
$1$ in $\lambda^{-}(v_{\beta})$ .
Thus $ v_{\beta} \in Z_i(\mathbf{w}_0) \; \mbox{if and only if} \, 
(\alpha_i, \lambda^-(v_{\beta}))>0.$

There is only one projective indecomposable module in $P_i(Q)$,
namely the projective cover of $S_i$, which is on the $i$-th
translation level of $\Gamma_Q$. Likewise, the left border
chamber of $\mathcal{WD}(\mathbf{w}_0)$ on level $j, \; j \in I$,
has the weight  $\omega_j$. Hence the only vertex on the left border
of $G^0(\mathbf{w}_0)$ verifying $(\alpha_i, \lambda^-(v_{\beta}))>0$,
is the one on the $i^{th}$ level. 

It remains to verify Proposition 5.4 for non-projective vertices
of $\Gamma_Q$. Lemma 5.5 gives us, for all $i,k$ that
$ (\alpha_k, \rho_i)=-(\alpha_i, c \rho_k).$
If $\beta=\underset{k=1}{\overset{n}{\sum}} m_k \alpha_k$, then
$$ \begin{array}{rl}
(\beta, \rho_i)& =
\underset{k=1}{\overset{n}{\sum}} m_k (\alpha_k, \rho_i) \\
&= -\underset{k=1}{\overset{n}{\sum}} m_k (\alpha_i, c \rho_k) \\
&= (\alpha_i, c (-\underset{k=1}{\overset{n}{\sum}} m_k \rho_k)) \\
&= (\alpha_i, c \lambda^+(v_{\beta}))\\
&= (\alpha_i, \lambda^- (v_{\beta})). 
\end{array}$$

We see that $(\beta, \rho_i)>0$ if and only if
$(\alpha_i, \lambda^-(v_{\beta}))>0$ so the 
first statement
of the proposition is verified. $\Box$

\section*{6. String cone inequalities and Lusztig moves}

\hspace*{5mm}
A vertex occurring in  a GP-path of type $i$ must belong to
a chamber of $\mathcal{Z}_i(\mathbf{w}_0)$. We shall
denote by $Y_i(\mathbf{w}_0)$ the set of these vertices. 
We have already accounted for the set 
$Z_i(\mathbf{w}_0)$ which consists of right-most vertices
of the chambers in $\mathcal{Z}_i(\mathbf{w}_0)$. It remains
to compute $Y_i(\mathbf{w}_0) \backslash Z_i(\mathbf{w}_0)$.
As we shall now see, these vertices lie on the border path $\delta_i$.

\vfill \eject

\textbf{Proposition 6.1}

\emph{
Fix $i \in I$. Then the subgraph with vertices
$Y_i(\mathbf{w}_0)$
is given inside $G(\mathbf{w}_0,i)$, according to the neighbourhood
of $i$ inside $Q$, 
by Figures 6.1 below  (with the relative position of $Z_i(\mathbf{w}_0)$
inside $Y_i(\mathbf{w}_0)$ delimited by segmented lines) :}

\[
\begin{xy}
(40,75)*{a)\stackrel{i-1}{\cdot} \longrightarrow
           \stackrel{i}{\cdot} \longleftarrow 
           \stackrel{i+1}{\cdot}};
(30,115)*{};
(35,115)*{} **\dir{-};
(40,120)*{\cdot} **\dir{-} ?(0.5)*@{>};
(40,120)*{\cdot};
(45,125)*{} **\dir{-} ?>*@{>};
(45,115)*{};
(40,120)*{\cdot} **\dir{-} ?(0.5)*@{>};
(40,120)*{\cdot};
(35,125)*{} **\dir{-};
(30,125)*{} **\dir{-} ?(0.5)*@{>};
(25,115)*{l_{i+1}};
(25,125)*{l_{i}};
(40,126)*{v_{\alpha_i}};
(40,123)*{};
(43,120)*{} **\dir{-};
(40,117)*{} **\dir{-};
(37,120)*{} **\dir{-};
(40,123)*{} **\dir{-};
(50,120)*{Z_i(\mathbf{w}_0)};
(120,75)*{b)\stackrel{i-1}{\cdot} \longleftarrow
           \stackrel{i}{\cdot} \longrightarrow 
           \stackrel{i+1}{\cdot}};
(80,115)*{l_{i+1}};
(80,125)*{l_{i}};
(85,115)*{};
(90,115)*{} **\dir{-} ?(0.5)*@{>};
(90,115)*{};
(95,110)*{\cdot} **\dir{-};
(105,100)*{\cdot} **\dir{-} ?(0.5)*@{>};
(105,100)*{\cdot};
(115,90)*{\cdot} **\dir{-} ?(0.5)*@{>};
(115,90)*{\cdot};
(120,85)*{} **\dir{-};
(130,85)*{} **\dir{-} ?(0.5)*@{>};
(130,85)*{};
(135,90)*{\cdot} **\dir{-};
(135,90)*{\cdot};
(145,100)*{\cdot} **\dir{-} ?(0.5)*@{>};
(145,100)*{\cdot};
(155,110)*{\cdot} **\dir{-} ?(0.5)*@{>};
(155,110)*{\cdot};
(160,115)*{} **\dir{-}; 
(165,115)*{} **\dir{-} ?(0.5)*@{>};
(155,110)*{\cdot};
(160,105)*{} **\dir{-};
(165,105)*{} **\dir{-} ?(0.5)*@{<};
(170,115)*{r_{i+1}};
(170,105)*{r_{i}};
(155,110)*{\cdot};
(145,120)*{\cdot} **\dir{-} ?(0.5)*@{>};
(145,120)*{\cdot};
(135,130)*{\cdot} **\dir{-} ?(0.5)*@{>};
(135,130)*{\cdot};
(125,140)*{\cdot} **\dir{-} ?(0.5)*@{>};
(125,140)*{\cdot};
(120,145)*{} **\dir{-};
(110,145)*{} **\dir{-} ?(0.5)*@{>};
(110,145)*{};
(105,140)*{\cdot} **\dir{-};
(95,130)*{\cdot} **\dir{-} ?(0.5)*@{>};
(95,130)*{\cdot};
(90,125)*{} **\dir{-};
(90,125)*{};
(85,125)*{} **\dir{-} ?(0.5)*@{>};
(90,105)*{};
(95,110)*{\cdot} **\dir{-};
(105,120)*{\cdot} **\dir{-};
(115,130)*{\cdot} **\dir{-};
(125,140)*{\cdot} **\dir{-};
(130,145)*{} **\dir{-}; ?>*@{>};
(100,95)*{};
(105,100)*{\cdot} **\dir{-};
(115,110)*{\cdot} **\dir{-};
(125,120)*{\cdot} **\dir{-};
(135,130)*{\cdot} **\dir{-};
(140,135)*{} **\dir{-}; ?>*@{>};
(110,85)*{};
(115,90)*{\cdot} **\dir{-};
(125,100)*{\cdot} **\dir{-};
(135,110)*{\cdot} **\dir{-};
(145,120)*{\cdot} **\dir{-};
(150,125)*{} **\dir{-}; ?>*@{>};
(140,85)*{};
(135,90)*{\cdot} **\dir{-};
(125,100)*{\cdot} **\dir{-};
(115,110)*{\cdot} **\dir{-};
(105,120)*{\cdot} **\dir{-};
(95,130)*{\cdot} **\dir{-};
(90,135)*{} **\dir{-} ?>*@{>};
(150,95)*{};
(145,100)*{\cdot} **\dir{-};
(135,110)*{\cdot} **\dir{-};
(130,115)*{Z_i(\mathbf{w}_0)};
(125,120)*{\cdot};
(115,130)*{\cdot} **\dir{-};
(105,140)*{\cdot} **\dir{-};
(100,145)*{} **\dir{-} ?>*@{>};
(103,120)*{};
(135,88)*{} **\dir{--};
(157,110)*{} **\dir{--};
(125,142)*{} **\dir{--};
(103,120)*{} **\dir{--};
(156,116)*{v_{\alpha_i}};
\end{xy}
\]

\[
\begin{xy}
(40,5)*{c)\stackrel{i-1}{\cdot} \longleftarrow
           \stackrel{i}{\cdot} \longleftarrow 
           \stackrel{i+1}{\cdot}};
(30,25)*{};
(45,25)*{} **\dir{-} ?(0.5)*@{>};
(45,25)*{};
(50,30)*{\cdot} **\dir{-};
(60,40)*{\cdot} **\dir{-} ?(0.5)*@{>};
(60,40)*{\cdot};
(70,50)*{\cdot} **\dir{-} ?(0.5)*@{>};
(70,50)*{\cdot};
(75,55)*{} **\dir{-} ?>*@{>};
(75,45)*{};
(70,50)*{\cdot} **\dir{-} ?(0.5)*@{>};
(70,50)*{\cdot};
(65,55)*{} **\dir{-};
(55,55)*{} **\dir{-} ?(0.5)*@{>};
(55,55)*{};
(50,50)*{\cdot} **\dir{-};
(40,40)*{\cdot} **\dir{-} ?(0.5)*@{>};
(40,40)*{\cdot};
(35,35)*{} **\dir{-};
(30,35)*{} **\dir{-} ?(0.5)*@{>};
(55,25)*{};
(50,30)*{\cdot} **\dir{-};
(40,40)*{\cdot} **\dir{-};
(35,45)*{} **\dir{-} ?>*@{>};
(65,35)*{};
(60,40)*{\cdot} **\dir{-};
(50,50)*{\cdot} **\dir{-};
(45,55)*{} **\dir{-} ?>*@{>};
(47,30)*{};
(50,27)*{} **\dir{--};
(73,50)*{} **\dir{--};
(70,53)*{} **\dir{--};
(47,30)*{} **\dir{--};
(25,25)*{l_{i+1}};
(25,35)*{l_{i}};
(77,50)*{v_{\alpha_i}};
(62,31)*{Z_i(\mathbf{w}_0)};
(120,5)*{d)\stackrel{i-1}{\cdot} \longrightarrow
           \stackrel{i}{\cdot} \longrightarrow 
           \stackrel{i+1}{\cdot}};
(95,45)*{};
(100,45)*{} **\dir{-} ?(0.5)*@{>};
(100,45)*{};
(105,40)*{\cdot} **\dir{-};
(115,30)*{\cdot} **\dir{-} ?(0.5)*@{>};
(115,30)*{\cdot};
(125,20)*{\cdot} **\dir{-} ?(0.5)*@{>};
(125,20)*{\cdot};
(130,15)*{} **\dir{-};
(140,15)*{} **\dir{-} ?(0.5)*@{>};
(140,15)*{};
(145,20)*{\cdot} **\dir{-};
(150,25)*{} **\dir{-};
(155,25)*{} **\dir{-} ?(0.5)*@{>};
(155,15)*{};
(150,15)*{} **\dir{-} ?(0.5)*@{>};
(150,15)*{};
(145,20)*{\cdot} **\dir{-};
(135,30)*{\cdot} **\dir{-} ?(0.5)*@{>};
(135,30)*{\cdot};
(125,40)*{\cdot} **\dir{-} ?(0.5)*@{>};
(125,40)*{\cdot};
(115,50)*{\cdot} **\dir{-} ?(0.5)*@{>};
(115,50)*{\cdot};
(110,55)*{} **\dir{-};
(95,55)*{} **\dir{-} ?(0.5)*@{>};
(100,35)*{};
(105,40)*{\cdot} **\dir{-};
(115,50)*{\cdot} **\dir{-};
(120,55)*{} **\dir{-} ?>*@{>};
(110,25)*{};
(115,30)*{\cdot} **\dir{-};
(125,40)*{\cdot} **\dir{-};
(130,45)*{} **\dir{-} ?>*@{>};
(120,15)*{};
(125,20)*{\cdot} **\dir{-};
(135,30)*{\cdot} **\dir{-};
(140,35)*{} **\dir{-} ?>*@{>};
(115,53)*{};
(112,50)*{} **\dir{--};
(145,17)*{} **\dir{--};
(148,20)*{} **\dir{--};
(115,53)*{} **\dir{--};
(152,20)*{v_{\alpha_i}};
(137,38)*{Z_i(\mathbf{w}_0)};
(90,45)*{l_{i+1}};
(90,55)*{l_{i}};
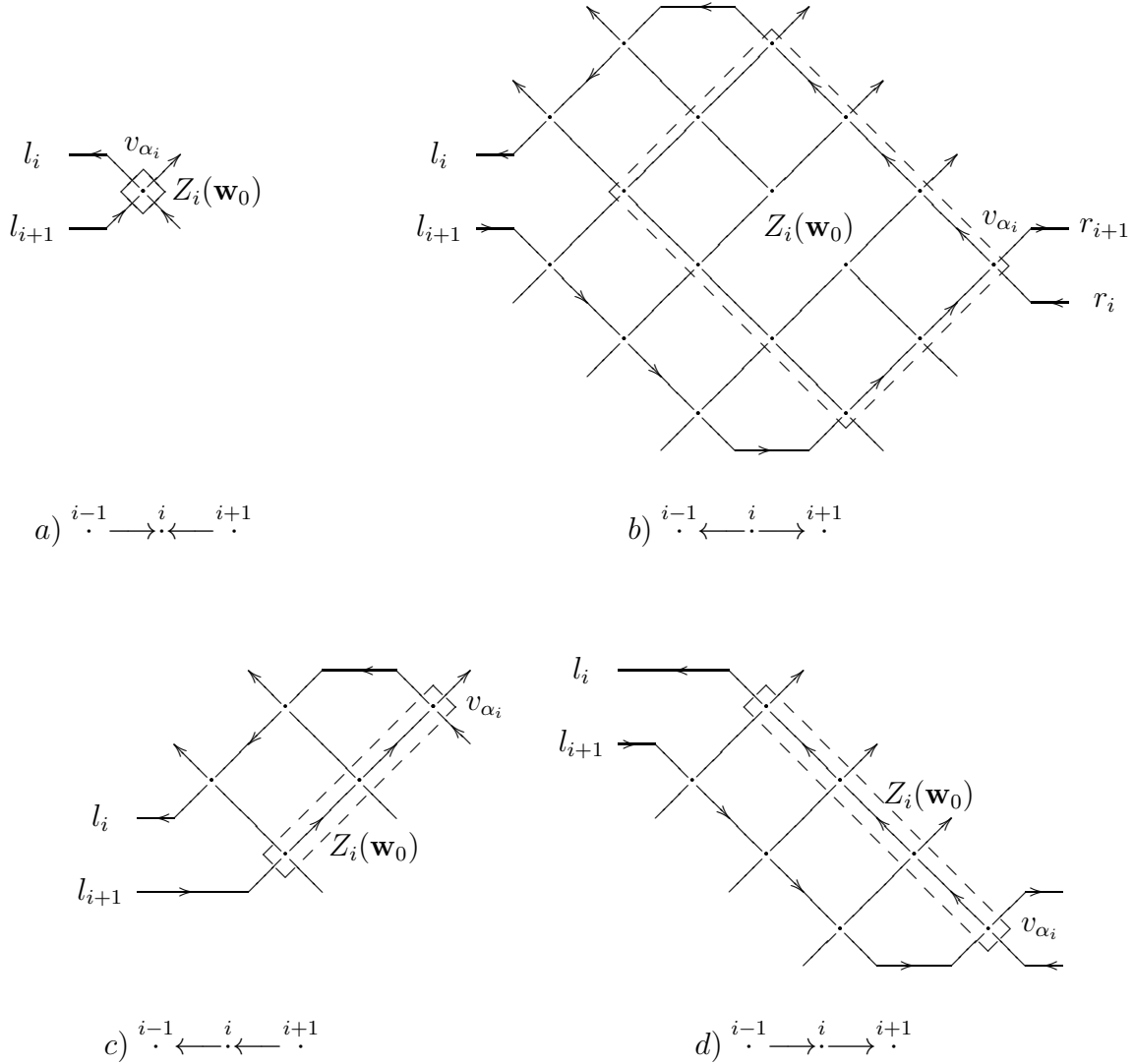
(90,-5)*{\mbox{Figures 6.1 : structure of $Y_i(\mathbf{w}_0)$ 
inside $G(\mathbf{w}_0,i)$.}};
\end{xy}
\]

\textbf{Proof :}

The key to the proof is the combinatorics of $H_i(Q)$ 
(and therefore of $H_i(\mathbf{w}_0)$) . Let
$(j_1, j_2 \ldots j_i  \mid j_{i+1} \ldots j_{n+1})$ be the
$i$-th segmented expression of $c$. 
The vertex of $H_i(Q)$ at position $(k,l)$ in Figure 5.1,
is mapped by $\Psi$ onto the intersection of the 
pseudo-lines $L_{j_k}$ and $L_{j_l}$, where 
$j_k \leq i$ and $j_l \geq i+1$.
The structure of $Z_i(\mathbf{w}_0)$ inside
$G(\mathbf{w}_0,i)$ follows then from 
Proposition 5.4. It remains to understand the position
of the border path $\delta_i$ with respect to
$Z_i(\mathbf{w}_0)$.

We shall call  the set of
vertices $v_{\beta} \in Z_i(\mathbf{w}_0)$ such that
either $\tau_{wd} (v_{\beta})$ is not defined, 
or  $\tau_{wd} (v_{\beta}) \not \in Z_i(\mathbf{w}_0)$,
the left border of $Z_i(\mathbf{w}_0)$.
We shall call the set of vertices $v_{\beta} \in Z_i(\mathbf{w}_0)$ 
which are not the image under $\tau_{wd}$ of another vertex
of $Z_i(\mathbf{w}_0)$,
the right border of $Z_i(\mathbf{w}_0)$. 
We shall see the limiting path $\delta_i$ consists
of vertices either lying on the right border
of $Z_i(\mathbf{w}_0)$ or of vertices
obtained by applying the translation $\tau_{wd}$
on a vertex lying on the left border of $Z_i(\mathbf{w}_0)$.
This ensures that vertices of $Y_i(\mathbf{w}_0)$ which
do not belong to $Z_i(\mathbf{w}_0)$, lie on $\delta_i$.

Let us consider the part $\delta_i^{<}$ of $\delta_i$ 
lying on the pseudoline $L_i$.

\underline{Case 1 : } One has 
$\stackrel{i-1}{\cdot} \longrightarrow \stackrel{i}{\cdot}$ inside $Q$.

The $i$-segmented expression of $c$ is such that
 $j_i=i$. In view of Proposition 5.2 b), $L_i$
exits $Z_i(\mathbf{w}_0)$ on the leftmost
vertex $v_{\beta_{min}}$ of $H_i(\mathbf{w}_0)$.
Now
$\lbrack \beta_{min} \rbrack$ is a projective
isoclass of $\mbox{mod} \, \Cbold Q$, hence 
$v_{\beta_{min}}$ is on the left border of 
$G^{\circ}(\mathbf{w}_0)$. The pseudoline $L_i$ goes
directly to $l_i$ after leaving
this vertex. The segment $\delta_i^{<}$ contains only
one vertex of $Z_i(\mathbf{w}_0)$, namely $v_{\beta_{min}}$,
which lies on the right border of this set. 

\underline{Case 2 : } One has 
$\stackrel{i-1}{\cdot} \longleftarrow \stackrel{i}{\cdot}$ 
inside $Q$.

This time  $j_1=i$, so $L_i$ leaves $Z_i(\mathbf{w}_0)$ 
on level 1. Recall $j_{i+1}, \ldots j_{n+1}$
is a permutation of the interval $i+1, i+2, \ldots ,n+1$,
so $Z_i(\mathbf{w}_0)$ contains all vertices 
$L_i \cap L_k, \; k  \geq i+1$, these vertices
being on the right border of $Z_i(\mathbf{w}_0)$.
 
Let us fix $k \leq i-1$, and apply
$\tau_{wd}$
on the leftmost vertex $L_{j_k} \cap L_{j_{n+1}}$
of $Z_i(\mathbf{w}_0)$ on level $k$. The result 
is, according to the $i$-segmented expression of $c$, the vertex
$L_{j_{k+1}} \cap L_{j_1}$. Yet $j_1=i$, and
$j_2, j_3, \ldots j_i$ is a permutation
of $1, 2, \ldots , i-1$. We see all
vertices $L_i \cap L_{k}$ with $k \leq i-1$ lie
on $\delta_i^{<}$. We 
have accounted for all vertices lying on $L_i$, 
so that after leaving $Y_i(\mathbf{w}_0)$,
$L_i$ goes to $l_i$.
  
The case of  the part $\delta_i^{>}$ of $\delta_i$ 
lying on strand $L_{i+1}$
is totally symmetric, so we won't give all details. 
The two cases to consider, correspond to
the orientations of the arrow between $i$ and $i+1$.
If one has 
$\stackrel{i}{\cdot} \longleftarrow \stackrel{i+1}{\cdot}$
inside $Q$, $j_{n+1}=i+1$ so $L_{i+1}$ enters
$Z_i(\mathbf{w}_0)$ on the first (and 
projective) vertex $v_{\beta_{min}}$,
directly from the border vertex $l_{i+1}$.
In the other case,
$\stackrel{i}{\cdot} \longrightarrow \stackrel{i+1}{\cdot}$
inside $Q$, $L_{i+1}$ enters 
$Z_i(\mathbf{w}_0)$ on level $n$, and 
we can account for all 
vertices $L_k \cap L_{i+1}$.
Either they lie on the right border of $Z_i(\mathbf{w}_0)$ if
$k \leq i$, or are translates of vertices at levels $k=i+2, \ldots n$ 
of the left border of that set. $\Box$
\medskip

\textbf{Proof of the Main Theorem 2.4}

We shall proceed in two steps : first we establish
a one-to-one correspondence \linebreak[4]
$\pi \mapsto A(\pi)$ 
between GP-paths of type $i$ and antichains of 
$P_i(Q)$, then fixing
a GP-path $\pi$, we show that the vector
$\mathbf{k}_\pi$ has the 
same coefficients as the 
Lusztig move $\mathbf{l}_{A(\pi)}$.

A GP-path $\pi$ of type $i$, consists by the oriented graph structure 
of $Y_i(\mathbf{w}_0)$ detailed in 
Figure 6.1, of three parts : A segment of $\delta_i^{>}$ from $l_{i+1}$ up to entry 
into $Z_i(\mathbf{w}_0)$, a path $\pi_ {grid}$ with vertices belonging to 
$Z_i(\mathbf{w}_0)$ up to an exit vertex, and finally return to 
$l_i$ on a segment of $\delta_i^{<}$. 
The first and last segments are uniquely defined by the first and
last vertices of $\pi_{grid}$, thus the segment $\pi_{grid}$
uniquely defines the whole path $\pi$. 

Let us use the coordinate system given in Proposition 5.2
for $P_i(Q)$ which lies inside $H_i(Q)$. The order structure
of $P_i(Q)$ is given by $(k_1, \, l_1) \leq (k_2, \, l_2)$
if and only if \linebreak[4] $k_1 \geq k_2 \; \mbox{and} \; l_1 \geq l_2$.
All the vertices of $Z_i(\mathbf{w}_0)$ are by
Propositions 5.2 and 5.4 of the same type :
a crossing of a line $L_{j_k}$ with $j_k \leq i$
which points backwards and upwards, 
with a line $L_{j_l}$ with $j_l \geq i+1$
which points forwards and upwards. The path
$\pi_{grid}$ is therefore a staircase path.
As such it is uniquely defined by the
"extremities" of its steps, vertices which are
a turning point from a forwards line
into a backwards line.
The inverse image under $\Psi$ of  these vertices
inside $P_i(Q)$ 
form a set $\{ (k_1, \, l_1), \, (k_2, \, l_2), \, \ldots ,\, 
(k_r,l_r) \}$ whose elements have coordinates which are pairwise disjoint : 
for any $j \neq j'$, 
one has $k_j \neq k_{j'}$ and $l_j \neq l_{j'}$.
By the nature of the order structure of $P_i(Q)$,
this set is an antichain that we shall denote 
by $A(\pi)$.

Conversely, an antichain of $P_i(Q)$ is a
set of vertices $\{ (k_1, \, l_1), \, (k_2, \, l_2), \, \ldots ,\, 
(k_r,l_r) \}$ that have coordinates which are pairwise disjoint.
The images of these vertices under $\Psi$
form the extremities of a uniquely defined
staircase path $\pi_{grid}$ 
of $Z_i(\mathbf{w}_0)$. By adding the appropriate
segments of the limiting path $\delta_i$,
we get a path $\pi(A)$ starting from $l_{i+1}$ and
returning to $l_i$. 
The nature of crossings inside $Z_i(\mathbf{w}_0)$
and Lemma 4.2 exclude the existence of forbidden crossings.
The path $\pi(A)$ is therefore a GP-path of type $i$.

The mapping $A \mapsto \pi(A)$ is the
inverse of  $\pi \mapsto A(\pi)$, and
vice-versa, so we have a one-to-one correspondence
$\psi$ between GP-paths of type $i$ and
antichains of $P_i(Q)$.

\[
\begin{xy}
(15,45)*{};
(30,60)*{\bullet} **\dir{.};
(20,40)*{};
(35,55)*{} **\dir{.};
(25,35)*{};
(35,45)*{\bullet} **\dir{.};
(40,50)*{} **\dir{.};
(30,30)*{\bullet};
(45,45)*{} **\dir{.};
(35,25)*{};
(50,40)*{} **\dir{.};
(15,45)*{};
(30,30)*{\bullet} **\dir{.};
(35,25)*{} **\dir{.};
(20,50)*{};
(40,30)*{} **\dir{.};
(25,55)*{};
(35,45)*{\bullet} **\dir{.};
(45,35)*{} **\dir{.};
(30,60)*{\bullet};
(50,40)*{} **\dir{.};
(75,42)*{};
(78,42)*{} **\dir{-};
(95,25)*{\circ} **\dir{-} ?(0.5)*@{>};
(95,25)*{\circ};
(100,30)*{\bullet} **\dir{-} ?(0.5)*@{>};
(100,30)*{\bullet};
(95,35)*{\circ} **\dir{-} ?(0.5)*@{>};
(95,35)*{\circ};
(105,45)*{\bullet} **\dir{-} ?(0.5)*@{>};
(105,45)*{\bullet};
(95,55)*{\circ} **\dir{-} ?(0.5)*@{>};
(95,55)*{\circ};
(100,60)*{\bullet} **\dir{-} ?(0.5)*@{>};
(100,60)*{\bullet};
(95,65)*{} **\dir{-};
(78,48)*{} **\dir{-} ?(0.5)*@{>};
(78,48)*{};
(75,48)*{} **\dir{-};
(70,42)*{L_{i+1}};
(70,48)*{L_i};
(85,45)*{};
(95,35)*{\circ} **\dir{.};
(100,30)*{\bullet};
(105,25)*{} **\dir{.};
(90,50)*{};
(110,30)*{} **\dir{.};
(105,45)*{\bullet};
(115,35)*{} **\dir{.};
(100,60)*{};
(120,40)*{} **\dir{.};
(95,25)*{\circ};
(100,20)*{} **\dir{--};
(105,25)*{} **\dir{--};
(106,20)*{L_{i+1}};
(85,45)*{};
(95,55)*{\circ} **\dir{.};
(90,40)*{};
(105,55)*{} **\dir{.};
(105,45)*{\bullet};
(110,50)*{} **\dir{.};
(100,30)*{\bullet};
(115,45)*{} **\dir{.};
(105,25)*{};
(120,40)*{} **\dir{.};
(104,30)*{+};
(109,45)*{+};
(104,60)*{+};
(45,55)*{P_i(Q)};
(114,55)*{Z_i(\mathbf{w}_0)};
(140,70)*{\bullet \; \mbox{positive contribution}};
(140,60)*{\circ \; \mbox{negative contribution}};
(30,10)*{\mbox{Antichain} \; A};
(90,10)*{\mbox{GP-path} \; \pi(A)};
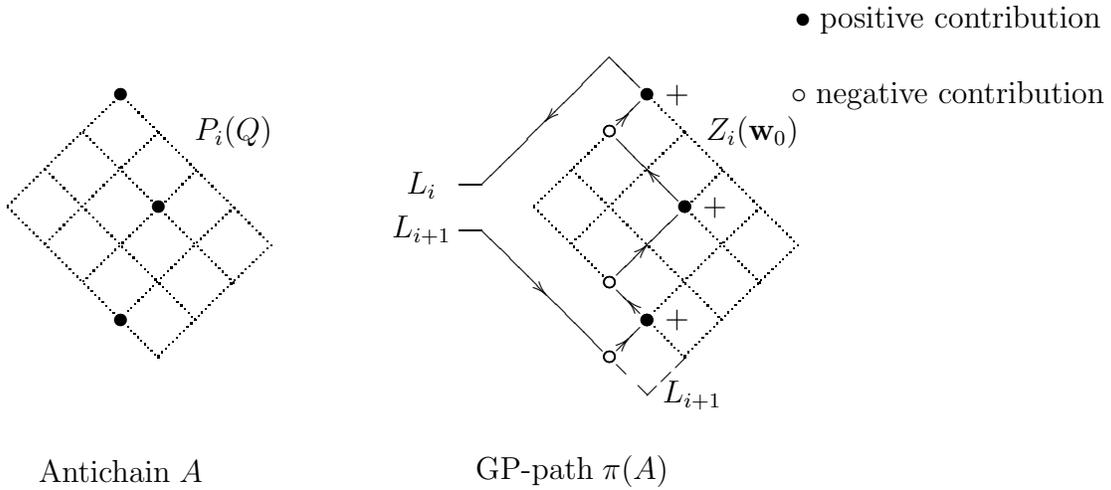
(90,0)*{\mbox{Figure 6.2 : Antichain of $P_i(Q)$ and 
corresponding GP-path of type $i$}};
\end{xy}
\]

Let us fix now a GP-path, and study its contributing
vertices. Recall a vertex $v_{\beta_h}$  contributes
a term $+t_h$ in $\mathbf{k}_{\pi} \cdot \tbold$ if $\pi$
changes strands $L_j \longrightarrow L_k$ at 
$v_{\beta_h}$ with a decrease from $j$ to $k$, and 
a term $-t_h$ 
if a change of strands $L_j \longrightarrow L_k$
occurs  at $v_{\beta_h}$
with an increase from $j$ to $k$.   
If $\pi$ stays on the same strand while passing
through $v_{\beta_h}$, no contribution occurs.

In view of Figure 6.2, a change of strands occurs
at an extremal vertex $v_{\beta_h}$ of $\pi_{grid}$.
At this vertex, 
the incoming strand $L_{j_k}$ verifies $j_k \geq i+1$,
while the outgoing strand $L_{j_l}$ verifies $j_l \leq i$,
We see a decrease of indices 
occurs in the passage $L_{j_k} \longrightarrow L_{j_l}$,
so $v_{\beta_h}$ is a positively contributing vertex.

Let us denote by $J_{\pi}$ the order ideal defined
by the extremal vertices of $\pi_{grid}$, inside $Z_i(\mathbf{w}_0)$ 
endowed with the poset structure induced by that $P_i(Q)$. 
Consider a minimal
element $v_{\beta_{h'}}$ of $Z_i(\mathbf{w}_0) \backslash J_{\pi}$,
and apply the translation $\tau_{wd}$. We get a vertex $\tau_{wd}(v_{\beta_{h'}})$ of $\pi$ 
with the following possible cases :
\begin{itemize}
\item[i)] $\tau_{wd}(v_{\beta_{h'}})$ is in $Z_i(\mathbf{w}_0)$, in which
case it is adjacent to three chambers of $\mathcal{Z}_i(\mathbf{w}_0)$.
A passage from a backward oriented line $L_{j_k}$ 
to a forward oriented line $L_{j_l}$ occurs, with 
an increase of indices.
\item[ii)] $\tau_{wd} (v_{\beta_ {h'}})$ is a vertex of $\delta_i^>$,
outside $Z_i(\mathbf{w}_0)$.  The strand $L_{i+1}$ is of minimal index among forwardly
oriented strands, so an increase of indices occurs by Proposition 6.1.
 
\item[iii)] $\tau_{wd}(v_{\beta_{h'}})$ is a vertex of 
$\delta_i^<$ outside $Z_i(\mathbf{w}_0)$. 
The strand $L_{i}$ is of maximal index among backwards
oriented strands, so again, by Proposition 6.1 an increase of indices 
occurs.
\end{itemize}
\[
\begin{xy}
(15,45)*{\cdot};
(25,55)*{\cdot} **\dir{.};
(35,65)*{\cdot} **\dir{.};
(25,35)*{\cdot};
(35,45)*{\circ} **\dir{.};
(45,55)*{\cdot} **\dir{-} ?(0.5)*@{>};
(35,25)*{\cdot};
(45,35)*{\cdot} **\dir{.};
(35,25)*{\cdot};
(25,35)*{\cdot} **\dir{.};
(15,45)*{\cdot} **\dir{.};
(45,35)*{\cdot};
(35,45)*{\circ} **\dir{-} ?(0.5)*@{>};
(35,45)*{\circ};
(25,55)*{\cdot} **\dir{.};
(45,55)*{\cdot};
(35,65)*{\cdot} **\dir{.};
(49,56)*{L_{j_l}};
(49,34)*{L_{j_k}};
(70,45)*{\cdot};
(80,35)*{\circ} **\dir{-} ?(0.5)*@{>};
(80,35)*{\circ};
(90,25)*{} **\dir{--};
(93,22)*{L_{i+1}};
(80,35)*{\circ};
(90,45)*{\cdot} **\dir{-} ?(0.5)*@{>};
(92,40)*{L_{j_l}};
(110,45)*{\cdot};
(100,35)*{\cdot} **\dir{.};
(90,45)*{\cdot} **\dir{.};
(80,55)*{\cdot} **\dir{.};
(90,65)*{\cdot} **\dir{.};
(95,50)*{Z_i(\mathbf{w}_0)};
(150,70)*{};
(140,60)*{\circ} **\dir{--};
(130,50)*{} **\dir{-} ?(0.5)*@{>};
(127,47)*{L_i};
(150,50)*{\cdot};
(140,60)*{\circ} **\dir{-} ?(0.5)*@{>};
(150,55)*{L_{j_k}};
(150,30)*{\cdot};
(140,40)*{\cdot} **\dir{.};
(150,50)*{\cdot} **\dir{.};
(160,60)*{\cdot} **\dir{.};
(170,50)*{\cdot} **\dir{.};
(155,45)*{Z_i(\mathbf{w}_0)};
(80,12)*{\mbox{Figure 6.3 \; : \; Negative
contributing vertices}};
\end{xy}
\]
We see that the translation $\tau_{wd}(v_{\beta_{h'}})$
results in all cases with a negative contribution
$-1$ to $\mathbf{k}_{\pi}$.

Finally, if $v_{\beta_h}$ is neither maximal inside
$J_{\pi}$, nor a translate of a minimal element
of $Z_i(\mathbf{w}_0) \backslash J_{\pi}$, $\pi$ stays on
the same strand while passing through it, so that 
no contribution occurs. 

Comparing this discussion with the definition of
the Lusztig move out of the \linebreak[4] antichain $A(\pi)$ we see
we get exactly the same positive and negative contributions
$\pm t_h$ in both cases. $\Box$

\textbf{Remarks}

i) The staircase paths used in the proof seem 
closely related to Le-diagrams defined by A. Postnikov 
in his study of positroids 
(compare Figure 6.2 with \cite{PostnikovPositroids} 
Figure 17.3). 

ii)  The neat structure of $Y_i(\mathbf{w}_0)$ described by Figure 6.1 
breaks down for reduced expressions $\mathbf{w}_0$ 
of type $A_n$ that are not adapted to a quiver.
Its subset $Z_i(\mathbf{w}_0)$ may include vertices which lie 
outside of $H_i(\mathbf{w}_0)$. Furthermore, the regular
grid structure of   $Z_i(\mathbf{w}_0)$ is destroyed
by the appearance of horizontal adjacencies 
between vertices of $Z_i(\mathbf{w}_0)$.

\section*{7. A $D_n$ example}

\hspace*{5mm} Let us consider the $D_4$ type quiver $Q$, given 
in the appendix of \cite{Reineke}  as satisfying condition $(L)$ :
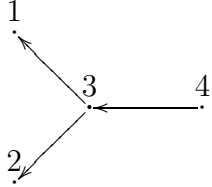
\[
\begin{xy}
(35,20)*{\cdot};
(20,20)*{\cdot} **\dir{-} ?>*@{>};
(20,20)*{\cdot};
(10,30)*{\cdot} **\dir{-} ?>*@{>};
(20,20)*{\cdot};
(10,10)*{\cdot} **\dir{-} ?>*@{>};
(10,33)*{1};
(10,13)*{2};
(20,23)*{3};
(35,23)*{4};
\end{xy}
\]

Let us verify
the conjecture stated in section 2. 
$\mathbf{w}_0=s_1 s_2 s_3 s_1 s_2 s_4 s_3 s_1 s_2 s_4 s_3 s_4$
is adapted to $Q$, with reflection ordering
$$\begin{array}{lll}
\beta_1=\alpha_1  &
\beta_5=\alpha_1+\alpha_3 & 
\beta_9=\alpha_2+\alpha_3+\alpha_4 \\
\beta_2=\alpha_2  & 
\beta_6=\alpha_1+\alpha_2+\alpha_3+\alpha_4 & 
\beta_{10}=\alpha_3 \\
\beta_3=\alpha_1+\alpha_2+\alpha_3  &
\beta_7=\alpha_1+\alpha_2+2\alpha_3+\alpha_4  & 
\beta_{11}=\alpha_3+\alpha_4 \\
\beta_4=\alpha_2+\alpha_3  &
\beta_8=\alpha_1+\alpha_3+\alpha_4  & 
\beta_{12}=\alpha_4
\end{array}$$

Lusztig parametrization is given by elements
$\tbold \in \Nbold^{12}$ with $t_i$ corresponding to 
$\lbrack \beta_i \rbrack$. We shall concisely denote
a vertex $\lbrack \beta \rbrack$ of $\Gamma_Q$ by
the indexes of the simple roots appearing. For instance
$134$ stands for $\lbrack \alpha_1+\alpha_3+\alpha_4 \rbrack$.
We shall denote the vertex
$\lbrack \alpha_1+\alpha_2+2 \alpha_3+\alpha_4 \rbrack$ 
by $12 \overline{3} 4$.  The Auslander-Reiten
quiver is given by
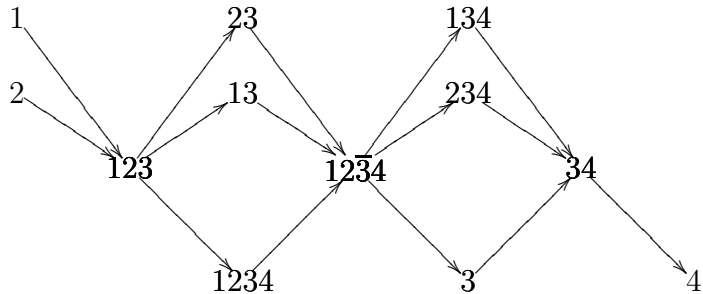
\[
\begin{xy}
(20,55)*{1};
(35,35)*{123} **\dir{-} ?>*@{>};
(35,35)*{123};
(50,55)*{23} **\dir{-} ?>*@{>};
(50,55)*{23};
(65,35)*{12\overline{3}4} **\dir{-} ?>*@{>};
(65,35)*{12\overline{3}4};
(80,55)*{134} **\dir{-} ?>*@{>};
(80,55)*{134};
(95,35)*{34} **\dir{-} ?>*@{>};
(20,45)*{2};
(35,35)*{123} **\dir{-} ?>*@{>};
(35,35)*{123};
(50,45)*{13} **\dir{-} ?>*@{>};
(50,45)*{13};
(65,35)*{12\overline{3}4} **\dir{-} ?>*@{>};
(65,35)*{12\overline{3}4};
(80,45)*{234} **\dir{-} ?>*@{>};
(80,45)*{234};
(95,35)*{34} **\dir{-} ?>*@{>};
(35,35)*{123};
(50,20)*{1234} **\dir{-} ?>*@{>};
(50,20)*{1234};
(65,35)*{12\overline{3}4} **\dir{-} ?>*@{>};
(65,35)*{12\overline{3}4};
(80,20)*{3} **\dir{-} ?>*@{>};
(80,20)*{3};
(95,35)*{34} **\dir{-} ?>*@{>};
(95,35)*{34};
(110,20)*{4} **\dir{-} ?>*@{>};
\end{xy}
\]
The $P_i(Q)$ sets are :
$$ \begin{array}{rl} 
P_1(Q)&=\{ 1 \} \\
P_2(Q)&=\{ 2 \} \\
P_3(Q)&=\{ 123, \, 23, \, 13, \,  12\overline{3}4, 3\} \\
P_4(Q)&=\{ 1234,\, 12\overline{3}4 , \, 134,
 \, 234, \, 34, \, 4 \}. 
\end{array}
$$  

There are respectively $1, 1, 6, 7$ antichains
of types $1, 2, 3, 4$.

Kashiwara's embedding corresponds to 
\cite{BerensteinZelevinskyInventiones} Theorem 5.10.
The computation of the set $K^{BZ}_{\mathbf{w}_0}$
requires  therefore, by 
\cite{BerensteinZelevinskyInventiones} Proposition 3.3 (iii),
the use of the reduced expression $\mathbf{w}_0^{op}$,
whose reflection ordering is reversed as compared
with that of $\mathbf{w}_0$ ($\alpha_4$ occurs first,
$\alpha_1$ occurs last), as well as a reversal 
in the numbering of coordinates. In our case,
the set of indices of $\mathbf{w}_0^{op}$ 
is $\ibold^{op}=(4,3,4,2,1,3,4,2,1,3,2,1)$.

Fix $i \in I$, and let $E(\omega_i)$ the corresponding
fundamental representation of $U_q(\gbold)$ of type $D_4$.
An $\mathbf{i^{op}}$-\textbf{trail} $\pi$ of type $i$ 
goes from $\omega_i$ to $w_0 s_i \omega_i$. It
is given by a set of coefficients $\mbold=(m_1, \, m_2, \ldots m_{12})$
such that the monomial 
$e_1^{m_1} e_2^{m_2} \ldots e_{12}^{m_{12}}$ induces
a non-zero mapping from the weight space 
$E(\omega_i)_{w_0 s_i \omega_i}$ to the
highest weight space $E(\omega_i)_{\omega_i}$. 
$\pi$ defines a sequence of weights
$\omega_i = \gamma_0, \, \gamma_1,  \, \ldots
\gamma_{12}=w_0 s_i \omega_i$ 
with 
$$ \gamma_k:=\gamma_0-\underset{l=1}{\overset{k}{\sum}}
m_l \alpha_{i_l}.$$

The trail $\pi$ defines  
a vector $\hbold_{\pi} \in \Zbold^{12}$
whose $k$-th coordinate   is 
$h_k:=(\dfrac{\gamma_{k-1}+\gamma_{k}}{2}, \alpha_{i_k})$,
$(k=1, \ldots, 12)$ (\cite{BerensteinZelevinskyInventiones} page 5, (2.2)).
Let us denote $\kbold_{\pi}:=(h_{12}, \, h_{11}, \, \ldots , \, h_1)$. Then, by
\cite{BerensteinZelevinskyInventiones} Theorem 3.10,
the set $K^{BZ}_{\mathbf{w}_0}$ of all vectors $\kbold_{\pi}$,
where  $\pi$ is any $\ibold^{op}$-trail, of any type $i$ defines 
$\mathcal{C}_{\mathbf{w}_0}$.

Let us consider the $\ibold^{op}$-trail $\pi$ of type 3 with
coefficients $(0,1,1,1,1,1,1,0,0,1,1,1)$. The list of weights $\gamma_k$ 
through which $\pi$ passes, as well as the coordinates $h_k$, are given
by the following table (weights are expressed by their coordinates with
respect to the basis of fundamental weights) :
$$ \begin{array}{rllcr}
k & \gamma_k & (\gamma_{k-1}+\gamma_k)/2 & \alpha_{i_k} & h_k \vspace*{2mm} \\
0 & \lbrack 0,0,1,0 \rbrack & & & \\
1 & \lbrack 0,0,1,0 \rbrack & \lbrack 0,0,1,0 \rbrack & \alpha_4 & 0 \\
2 & \lbrack 1,1,-1,1 \rbrack & \lbrack \tfrac{1}{2}, \tfrac{1}{2}, 0 , \tfrac{1}{2} \rbrack & \alpha_3 & 0 \\
3 & \lbrack 1,1,0,-1 \rbrack & \lbrack 1,1,-\tfrac{1}{2}, 0 \rbrack & \alpha_4 & 0 \\
4 & \lbrack 1,-1,1,-1 \rbrack & \lbrack 1,0,\tfrac{1}{2},-1 \rbrack & \alpha_2 & 0 \\
5 & \lbrack -1,-1,2,-1 \rbrack & \lbrack 0,-1,\tfrac{3}{2},-1 \rbrack & \alpha_1 & 0 \\
6 & \lbrack 0,0,0,0 \rbrack & \lbrack -\tfrac{1}{2},-\tfrac{1}{2},1,-\tfrac{1}{2} \rbrack & \alpha_3 & 1 \\
7 & \lbrack 0,0,1,-2 \rbrack & \lbrack 0,0,\tfrac{1}{2}, -1 \rbrack & \alpha_4 & -1 \\
8 & \lbrack 0,0,1,-2 \rbrack & \lbrack 0,0, 1, -2 \rbrack & \alpha_2 & 0 \\ 
9 & \lbrack 0,0,1,-2 \rbrack & \lbrack 0,0,1,-2 \rbrack & \alpha_1 & 0 \\
10 & \lbrack 1,1,-1,-1 \rbrack & \lbrack \tfrac{1}{2}, \tfrac{1}{2},0,-\tfrac{3}{2} \rbrack 
& \alpha_3 & 0 \\
11 & \lbrack 1,-1,0,-1 \rbrack & \lbrack 1,0,-\tfrac{1}{2}, -1 \rbrack & \alpha_2 & 0 \\
12 & \lbrack -1,-1,1,-1 \rbrack & \lbrack 0,-1,\tfrac{1}{2},-1 \rbrack & \alpha_1 & 0 
\end{array}$$
We get $\hbold_{\pi}=(0,0,0,0,0,1,-1,0,0,0,0,0)$, hence
$\kbold_{\pi}=(0,0,0,0,0,-1,1,0,0,0,0,0)$, which defines
the string cone inequality $t_7-t_6 \geq 0$. If we consider
the antichain $A=\{ 12\overline{3}4 \}$ of $P_3(Q)$, then
$V_A=\lbrack \beta_7 \rbrack, \; U_A=\lbrack \beta_6 \rbrack$,
so we get $\lbold_A=\kbold_{\pi}$.

The enumeration of all possible $\ibold^{op}$-trails in our case
(with the help of quagroup package \cite{Quagroup}) shows
there is a one-to-one correspondence
between antichains $A$ of $P_i(Q)$ and $\ibold^{op}$-trails
$\pi$ of type $i$, for each of the types $i=1, \ldots , 4$. 
For every antichain $A$ there is a unique trail $\pi$ with 
$\lbold_A=\kbold_{\pi}$.
 We give these correspondences in the table below.
The first column is 
the set of 
defining inequalities of $\mathcal{C}_{\mathbf{w}_0}$.
Next to each inequality is the antichain $A$ defining 
the corresponding Lusztig move $\lbold_A$, and the $\ibold^{op}$-trail leading
to the corresponding vector 
$\kbold_{\pi} \in K^{BZ}_{\mathbf{w}_0}$.
\medskip

\begin{tabular}{c|r|c|c|c}
Type & Inequality & Antichain 
& $\ibold^{op}$-trail position &
$\ibold^{op}$-trail coefficients $\mathbf{m}$ \\
\hline
1 & $t_1 \geq 0$ & 1 & 
$4 \, 3 \, 4 \, 2 \, \underline{1} \, \underline{3} \, 
\underline{4} \, \underline{2} \, 1 \, \underline{3} \, 2 \, 1$ &
(0,0,0,0,1,1,1,1,0,1,0,0) \\
\hline
2 & $t_2 \geq 0$ & 2 & 
$4 \, 3 \, 4 \, \underline{2} \, 1 \, \underline{3} \,
\underline{4} \, 2 \, \underline{1} \, \underline{3} \,2 \, 1$ &
(0,0,0,1,0,1,1,0,1,1,0,0) \\
\hline
3 & $t_3-t_1-t_2 \geq 0$ & 123 & 
$4 \, \underline{3} \, \underline{4} \, \underline{2} \, \underline{1} \,
\underline{3} \, \underline{4} \, \underline{2} \, \underline{1} \,
3 \, 2 \, 1 $ & (0,1,1,1,1,2,1,1,1,0,0,0) \\
 & $t_4-t_2 \geq 0$ & 23 & 
$ 4 \, \underline{3} \, \underline{4} \, \underline{2} \, \underline{1} \,
\underline{3} \, \underline{4} \, \underline{2} \, 
1 \, 3 \, 2 \, \underline{1}$ & (0,1,1,1,1,2,1,1,0,0,0,1) \\
 & $t_5-t_1 \geq 0$ & 13 & 
$4 \, \underline{3} \, \underline{4} \, \underline{2} \, \underline{1} \,
\underline{3} \, \underline{4} \, 2 \, \underline{1} \,
3 \, \underline{2} \, 1 $ & (0,1,1,1,1,2,1,0,1,0,1,0) \\
& $t_4+t_5-t_3 \geq 0$ & 13, 23 &
$4 \, \underline{3} \, \underline{4} \, \underline{2} \, \underline{1} \,
\underline{3} \, \underline{4} \, 2 \, 1 \,
3 \, \underline{2} \, \underline{1} $ & (0,1,1,1,1,2,1,0,0,0,1,1) \\
 & $t_7-t_6 \geq 0$ & $12\overline{3}4$ & 
$4 \, \underline{3} \, \underline{4} \, \underline{2} \, \underline{1} \,
\underline{3} \, \underline{4} \, 2 \, 1 \,
\underline{3} \, \underline{2} \, \underline{1} $ & (0,1,1,1,1,1,1,0,0,1,1,1) \\
& $t_{10} \geq 0$ & 3 & 
$4 \, \underline{3} \, 4 \, \underline{2} \, \underline{1} \,
\underline{3} \, \underline{4} \, 2 \, 1 \,
\underline{3} \, \underline{2} \, \underline{1} $ & (0,1,0,1,1,1,2,0,0,1,1,1) \\
\hline
4 & $t_6-t_3 \geq 0$ & 1234 & 
$\underline{4} \, \underline{3} \, 4 \, \underline{2} \, \underline{1} \, 
\underline{3} \, 4 \, 2 \, 1 \, 3 \, 2 \, 1$  &
(1,1,0,1,1,1,0,0,0,0,0,0) \\
  & $t_7-t_4-t_5 \geq 0$ & $12\overline{3}4$ & 
$\underline{4} \, \underline{3} \, 4 \, \underline{2} \, \underline{1} \, 3 \,
4 \, 2 \, 1 \, \underline{3} \, 2 \, 1$  & (1,1,0,1,1,0,0,0,0,1,0,0) \\
 & $t_8-t_5 \geq 0$ & 134 &
$\underline{4} \, \underline{3} \, 4 \, \underline{2} \, 1 \, 3 \,
4 \, 2 \, \underline{1} \, \underline{3} \, 2 \, 1$&
(1,1,0,1,0,0,0,0,1,1,0,0) \\
 & $t_9-t_4 \geq 0$ & 234 &
$\underline{4} \, \underline{3} \, 4 \, 2 \, \underline{1} \, 3 \,
4 \, \underline{2} \, 1 \, \underline{3} \, 2 \, 1$ &
(1,1,0,0,1,0,0,1,0,1,0,0) \\
 & $t_8+t_9-t_7 \geq 0$ & 134, \, 234 &  
$\underline{4} \, \underline{3} \, 4 \, 2 \, 1 \, 3 \,
4 \, \underline{2} \, \underline{1} \, \underline{3} \, 2 \, 1$ 
& (1,1,0,0,0,0,0,1,1,1,0,0) \\
 & $t_{11}-t_{10} \geq 0$ & 34 &
$\underline{4} \, 3 \, 4 \, 2 \, 1 \, \underline{3} \,
\, 4 \underline{2} \, \underline{1} \, \underline{3} \, 2 \, 1$ &
(1,0,0,0,0,1,0,1,1,1,0,0) \\
 & $t_{12} \geq 0$ & 4 & 
$4 \, 3 \, \underline{4} \, 2 \, 1 \, \underline{3} \, 
4 \, \underline{2} \, \underline{1} \, \underline{3} \, 2 \, 1$ &
(0,0,1,0,0,1,0,1,1,1,0,0)
\end{tabular}
\medskip

\textbf{Remark : } The trail $\pi$ of type 3 considered above in detail,
passes through the weight $\lbrack 0,0,0,0 \rbrack$. All other trails
of the table pass only through extremal weights. The existence of $\pi$
 shows that  condition 
$(L)$ is not strong enough to allow an analogue of 
\cite{BerensteinZelevinskyInventiones} Theorem 3.14. The subexpression corresponding to
$\pi$ is not a reduced word.

\begin{center}
{ \small S. Zelikson, LMNO,  UMR 6139 du CNRS  } \\ 
{ \small D\'{e}partement de Math\'{e}matiques} \\ 
{ \small Universit\'{e} de Caen B.P. 5186 } \\
{ \small 14032 Caen Cedex FRANCE} \\
{ \small E-mail : Shmuel.Zelikson@unicaen.fr}
\end{center}

\end{document}